\documentclass[10pt,journal, onecolumn]{IEEEtran}

\usepackage{color,latexsym,amsfonts,amssymb,bbm,comment,amsmath,cite, amsthm}
\usepackage{hyperref}
\usepackage{fullpage}
\usepackage{graphicx}
\usepackage{tikz}
\usetikzlibrary{decorations.pathreplacing}
\usepackage{mathtools}
\usepackage{ifthen}
\usepackage{babel}

\usepackage{adjustbox}%Usedtoconstrainimagestoamaximumsize
\usepackage{xcolor}%Allowcolorstobedefined
\usepackage{enumerate}%Neededformarkdownenumerationstowork
\usepackage{amsmath}%Equations
\usepackage{amssymb}%Equations
\usepackage{fancyvrb}%verbatimreplacementthatallowslatex
\usepackage{grffile}%extendsthefilenameprocessingofpackagegraphics
\usepackage{hyperref}
\usepackage{longtable}%longtablesupportrequiredbypandoc>1.10
\usepackage{booktabs}%tablesupportforpandoc>1.12.2
\usepackage[inline]{enumitem}%IRkernel/reprsupport(itusestheenumerate*environment)
\usepackage[normalem]{ulem}
\usepackage{float}
\usetikzlibrary{calc} 
\usetikzlibrary{arrows}
\pgfmathsetseed{14}

\newcommand{\coxSegment}[6]{

	\pgfmathsetmacro{\Xa}{#1}
	\pgfmathsetmacro{\Xb}{#2}
	\pgfmathsetmacro{\Xc}{#3}
	\pgfmathsetmacro{\Xd}{#4}

	\draw (\Xa - \factor*\Xc + \factor*\Xa,\Xb  - \factor*\Xd + \factor*\Xb) -- (\Xc + \factor*\Xc - \factor*\Xa, \Xd + \factor*\Xd - \factor*\Xb);
	\foreach \y in {1,2,...,#5}{
		\pgfmathsetmacro{\Xf}{random()}
		\pgfmathsetmacro{\Xg}{1- \Xf)}
		\pgfmathsetmacro{\Xh}{\Xa*\Xf + \Xc*\Xg}
		\pgfmathsetmacro{\Xi}{\Xb*\Xf+ \Xd*\Xg}
		\pgfmathparse{\Xh*\Xh + \Xi*\Xi > #6 ? 1: 0}
		\ifthenelse{\pgfmathresult>0}{\fill (\Xh, \Xi) circle (2pt);}{}
	}

}

\def\eq{\begin{equation}}
\def\en{\end{equation}}

\newtheorem{theorem}{Theorem}[section]
\newtheorem{lemma}[theorem]{Lemma}
\newtheorem{proposition}[theorem]{Proposition}
\theoremstyle{definition}
\newtheorem{definition}[theorem]{Definition}

\newcommand{\muv}{\mu_{\mathrm v}}

\newcommand{\Tc}{T_{\mathrm c}}
\newcommand{\lac}{\la_{\mathrm c}}
\newcommand{\rhoc}{\rho_{\mathrm c}}
\newcommand{\vc}{v_{\mathrm c}}
\newcommand{\vcc}{v^{\mathrm c}}
\newcommand{\rc}{r_{\mathrm c}}

\newcommand{\supp}{{\mathrm{supp}}}
\usepackage{color}

\newcommand{\vmin}{v_{\mathrm{min}}}
\newcommand{\vmax}{v_{\mathrm{max}}}
\newcommand{\la}{\lambda}
\newcommand{\dist}{\mathrm{dist}}
\newcommand{\one}{\mathbbmss{1}}
\def \k{\kappa}

\def\e{{\varepsilon}}
\def\eps{\varepsilon}

\newcommand{\R}{\mathbb R}
\newcommand{\Q}{\mathbb Q}
\newcommand{\E}{\mathbb E}
\def\S{\mathbb S}
\def\P{\mathbb P} 
\newcommand{\Z}{\mathbb Z}
\newcommand{\N}{\mathbb N}
\newcommand{\C}{\mathcal{C}}
\def\a{\alpha}
\def\d{{\mathrm d}}
\newcommand{\tra}{\mathcal{T}}

\begin{document}
\title{Connectivity in mobile device-to-device networks in urban environments}

%\affil[1]{Weierstrass Institute for Applied Analysis and Stochastics, Mohrenstra{\ss}e 39, 10117 Berlin, Germany}
%\affil[2]{Orange Labs, 44 Avenue de la R\'epublique, 92320 Ch\^atillon, France}

\author{
Elie~Cali,~\IEEEmembership{Orange,}
Alexander~Hinsen,~\IEEEmembership{WIAS,}
        Benedikt~Jahnel,~\IEEEmembership{TU Braunschweig, WIAS,}
        and Jean-Philippe~Wary,~\IEEEmembership{Orange}% <-this % stops a space
\thanks{Orange Labs, 44 Avenue de la République, 92320 Châtillon, France.}% <-this % stops a space
\thanks{Weierstrass Institute for Applied Analysis and Stochastics, Mohrenstraße 39, 10117 Berlin, Germany.}
\thanks{Technische Universit\"at Braunschweig, Institut f\"ur Mathematische Stochastik, Universit\"atsplatz 2, 38106 Braunschweig, Germany.}}% <-this % stops a space

\maketitle

\begin{abstract}
In this article we setup a dynamic device-to-device communication system where devices, given as a Poisson point process, move in an environment, given by a street system of random planar-tessellation type, via a random-waypoint model. Every device independently picks a target location on the street system using a general waypoint kernel, and travels to the target along the shortest path on the streets with an individual velocity. Then, any pair of devices becomes connected whenever they are on the same street in sufficiently close proximity, for a sufficiently long time. 
After presenting some general properties of the multi-parameter system, we focus on an analysis of the clustering behavior of the random connectivity graph. 
%We exhibit regimes in which the graph features an infinite connected component with probability zero or with positive probability.
%The presence, respectively absence, of percolation, can be seen as an indication for the ability, respectively disability, of the system to transmit messages over long distances. Let us mention that such an analysis can be useful in settings going far beyond the device-to-device use case, and can for example help to describe disseminations of communicable diseases in mobile populations.   
In our main results we isolate regimes for the almost-sure absence of percolation if, for example, the device intensity is too small, or the connectivity time is too large. 
%, the system features bad connectivity in the sense of absence of percolation. 
%Especially the case of large speeds is interesting since devices can potentially travel far while spending less time in close proximity to other devices. 
On the other hand, we exhibit parameter regimes of sufficiently large intensities of devices, under favorable choices of the other parameters, such that percolation is possible with positive probability. Most interestingly, we also show an in-and-out of percolation as the velocity increases.
The rigorous analysis of the system mainly rests on comparison arguments with simplified models via spatial coarse graining and thinning approaches. Here we also make contact to geostatistical percolation models with infinite-range dependencies. 

%Additionally, in this article, we give a first detailed description of the simulation code for the random connectivity graph. This simulator will be used to perform supporting numerical analysis for the system in bounded domains and the results of these simulations will be presented in the third report. 
\end{abstract}

%\input{Brainstorming}

%====================================================
%====================================================
%====================================================

\section{Introduction and Setting}\label{sec_intro}
The ever increasing demand for fast and reliable data exchange in communication systems poses great challenges but also offers opportunities for network operators around the globe. One very important aspect here is the exponentially increasing use of mobile devices such as cellphones or communicating cars. This is also reflected in the 5G specifications, where faster connections, higher throughput, and more capacity are envisioned via an enhanced mobile broadband, and ultra-reliable low-latency communications should enable the system to support time-crucial applications such as car-to-car communication. 

In this context, {\em device-to-device (D2D) communications} is considered one of the key concepts pervading a highly diverse set of use cases. On the one hand, D2D systems hold the potential to relief present day cellular networks from at least some of the system's pressure. On the other hand, D2D communications can provide for example faster and more robust connectivity. However, from an operator's perspective, D2D systems are much less controllable than traditional cellular network due to the dependence on the individual user behavior. This lack of control becomes even more severe, when the devices are {\em mobile}. Thus, in order to properly predict the performance of D2D systems, either with or without additional infrastructure, a detailed and comprehensive modeling and model analysis of such {\em mobile ad-hoc networks} (MANets) is imperative. Here, a natural ansatz is to incorporate the uncertainties of the system with the help of {\em probability theory} and more precisely {\em stochastic geometry}~\cite{baccelli2009stochastic1,baccelli2009stochastic2,haenggi2012stochastic,jahnel2020probabilistic}. The starting point is the modeling of {\em random initial locations} of the devices and their individual movements in cities. Based on this, the {\em transmission mechanism} between any pair of devices has to be represented with a reasonable level of detail, for example not neglecting that the devices need to spend some minimal amount of time in close proximity. Here also different operating systems (e.g., Android versus IOS systems) and their specific properties, as well as fluctuations in the transmission rates play an important role. Ultimately, then the {\em distribution of data} through the D2D system over time can be analyzed, leading to predictions for examples for the {\em connectivity, throughput,} or more generally, the {\em quality-of-service} of the network. 

Against the background of the current pandemic crisis of the SARS-CoV-2 {\em human-to-human virus}, let us highlight that such mobile D2D networks and the associated mathematical theory of stochastic geometry can also be employed (and partially is already in use) to model and analyze connectivity models in the context of the {\em spatial epidemiology of communicable diseases}. In fact, the approach and results presented in this article can be easily adjusted and interpreted with respect to human-to-human connectivity and the associated spreading of a contagious virus in a population, especially with the application of contact-tracing. The general connection mechanism outline in this manuscript features a relatively realistic setting for the human-to-human transmission of a contagious viruses, but needs appropriate parameter adjustments. In systems with more complex interaction, i.e. virus transmission, this model can serve as a lower bound to see if long-range communication is possible.

We start by describing our general model for a mobile multilayered D2D telecommunication network in detail. 
The first network layer, described in Section \ref{route_sec}, is given by the street system. In Section~\ref{dev_sec} we describe the system of initial device positions on the street system. In Section \ref{mob_sec} we introduce the paradigmatic mobility model for the devices given by the random waypoint model. Finally, in Section~\ref{con_sec} we introduce our notion of connectivity in the system. 

%%%%%%%%%%%%%%%%%%%%
%Route systems
%%%%%%%%%%%%%%%%%%%%
\subsection{Street systems via random segment processes}\label{route_sec}
\begin{figure}[!htpb]
\centering
 \input{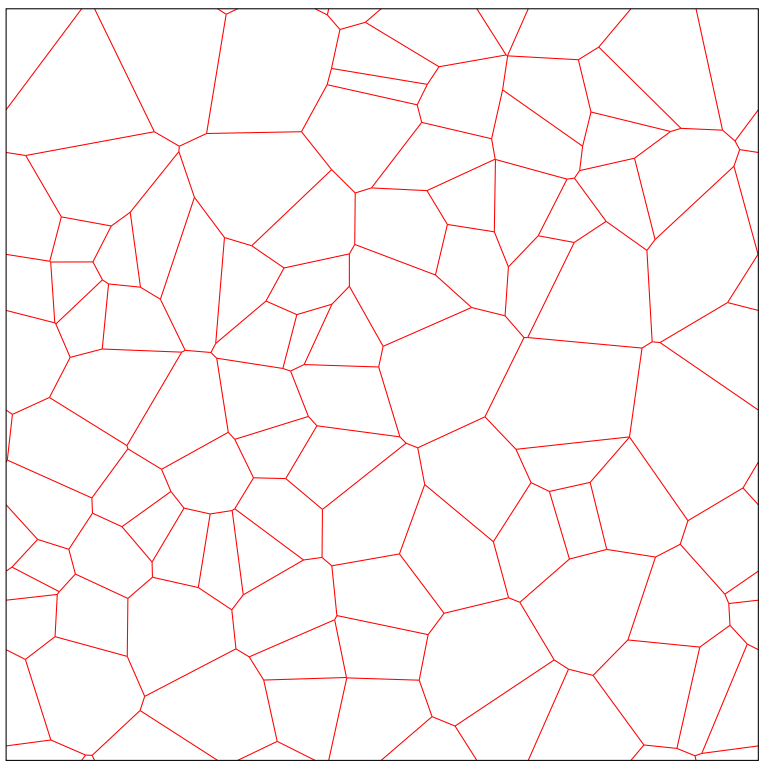}
\caption{Realization of a street system given by a Poisson--Voronoi tessellation.}
\label{Fig_Streets}
\end{figure}

The first model layer consists of a stationary random planar segment process $S\subset\R^2$, i.e., a simple point process of line segments of lengths in $(0,\infty)$ embedded in $\R^2$ (see~\cite{pp1}), that represents a {\em street system}. The most common choice for the street system model is the {\em Poisson--Voronoi tessellations (PVT)} based on statistical analysis of urban maps ~\cite{courtat2012promenade}. However, some of our results also cover for example the case of  {\em Poisson--Delaunay tessellations (PDT)}, which shares with the PVT the feature of translation-invariance and isotropy, i.e., the distribution of $S$ is invariant both under translations as well as rotations in $\R^2$, see Figure~\ref{Fig_Streets} for an illustration. 
We recall that the edges of a PVT are defined as the boundaries of the Voronoi cells of a stationary Poisson point process, where such a Poisson--Voronoi cell associated to the Poisson point $\Phi_i$ is a (almost-surely bounded) convex subset of $\R^2$ containing all points that are closer to $X_i$ than to any other Poisson point, see~\cite{okabe2000spatial} for more details. Another relevant choice for which some of our results are valid, is the {\em Manhattan grid} or {\em rectangular Poisson line tessellation}, where infinite horizontal and vertical lines are attached to two independent Poisson processes on the axis in $\R^2$. 
%Our tessellation processes can be parametrized for example by the intensity $\g>0$ of the street system, that is, $\g$ represents the expected length of streets per unit area~\cite{stoyanstochastic}. As a one-parameter description is not possible and from the perspective of an operator the street system is not changeable, we will try to suppress the street system wherever it is possible in the notation. 
We denote by $\mathcal{S}$ the law of the street system.% and indicate by $\gamma= \E\big[|S \cap [0,1]^2|\big]$ the expected length of street per unit volume.

%In order to apply and compare the results for different street systems, see Proposition~\ref{prop_scale} for a possible scaling relation.
% \red{NEW\\}
% While for PVT and PDT one parameter is enough to describe the whole law of the street system, we see $\gamma$ as a scaling parameter. 
% From the perspective of an operator, the law of the street system can not be changed. In order to streamline the notation, we are going to suppress the law of the street system and assume that for our street system $S$ we have $\gamma = 1$.
% For any street system $S^\gamma$ with $\gamma \neq 1$ we adjust our results by re-scaling the other parameters. The new street system $S^\gamma$ is obtained by the scaling $\mathcal{L}(d S^\gamma) = \mathcal{L}({\d \gamma S}$

% \red{END NEW}

It will become important that the underlying street system exhibits a certain spatial stochastic independence known as \emph{stabilization}, which has its roots in central limit theorems in stochastic geometry~\cite{stab1, stab2, stab3}. Loosely speaking, for our main results to hold, in two distant regions in space, the distribution of streets should not depend on one another. More precisely, we write $S\cap B$ to indicate the restriction of the street system $S$ to the area $B \subset \R^2$ and $|S\cap B|$ for the total lengths of the street system in $B$. Further, we write 
$$Q_n(x) = x + [-n/2, n/2]^2$$ 
for the square with side length $n\ge1$ centered at $x \in \R^2$ and put $Q_n = Q_n(o)$. We will assume throughout the manuscript that $S$ is normalized, i.e., $\E[|S\cap Q_1|]=1$. We write $\dist(\varphi, \psi) = \inf\{|x-y|:\, x\in\varphi, y\in\psi\}$ to denote the distance between sets $\varphi, \psi\subset\R^2$. Then, we reproduce the definition from \cite{wias2}.%\red{self-ref or more general?}

\begin{definition}
	\label{stabDef}
A stationary random segment process $S$ is called \emph{stabilizing}, if there exists a random field of \textit{stabilization radii} $R = \{R_x\}_{x\in\R^2}$ that is measurable with respect to $S$ and satisfies
\begin{enumerate}
 \item[(1)] $(S, R)$ are jointly stationary, 
 \item[(2)] $\lim_{n\uparrow\infty}\P\big(\sup_{y\in Q_n \cap\, \Q^2}R_y < n\big) = 1$,
 \item[(3)] for all $n \ge 1$, all non-negative bounded measurable functions $f$ and all finite $\varphi \subset \R^2$ with $\dist(x, \varphi \setminus \{x\}) > 3n$ for all $x \in \varphi$, the following random variables are independent
 $$\Big\{f(S\cap Q_n(x))\one\big\{\sup_{y \in Q_n(x) \cap \, \Q^2}R_y < n\big\}\Big\}_{x \in \varphi}.$$
\end{enumerate}
A stabilizing segment process $S$ is called \emph{exponentially stabilizing} if there exists $c>0$ such that $\P(\sup_{y\in Q_n \cap\, \Q^2}R_y \ge n)\le \exp(-cn)$ for all sufficiently large $n$.
\end{definition}
Note that PVT and PDT are even exponentially stabilizing, see for example~\cite{wias2}, while the Manhattan grid is not stabilizing. In the following, we will always assume that $\mathcal S$ is stabilizing. 

\medskip
Next, it will also become important that the underlying street system exhibits a sufficient amount of connectedness with high probability. We encode this in a notion of {\em asymptotically essentially connectedness} as follows. 
\begin{definition}
    \label{aecDef}
A stabilizing random segment process $S$ with stabilization radii $R$ is \emph{asymptotically essentially connected} if for all sufficiently large $n \ge 1$, whenever $\sup_{y\in Q_{2n}\cap\, \Q^2}R_y<n/2$ we have that 
\begin{enumerate}
            \item[(1)] $|S\cap Q_n|>0$ and 
            \item[(2)] $S\cap Q_n$ is contained in one of the connected components of $S\cap Q_{2n}$.
\end{enumerate}
\end{definition}
Again, street systems $S$ that are given as edges of a PVT or PDT, are asymptotically essentially connected, see for example~\cite[Example 3.1]{wias2} and in the following, we will always assume that $\mathcal S$ is asymptotically essentially connected.

%%%%%%%%%%%%%%%%%%%%
%Devices
%%%%%%%%%%%%%%%%%%%%
\subsection{Initial device positions via Cox point processes}\label{dev_sec}
\begin{figure}[!htpb]
\centering
 \input{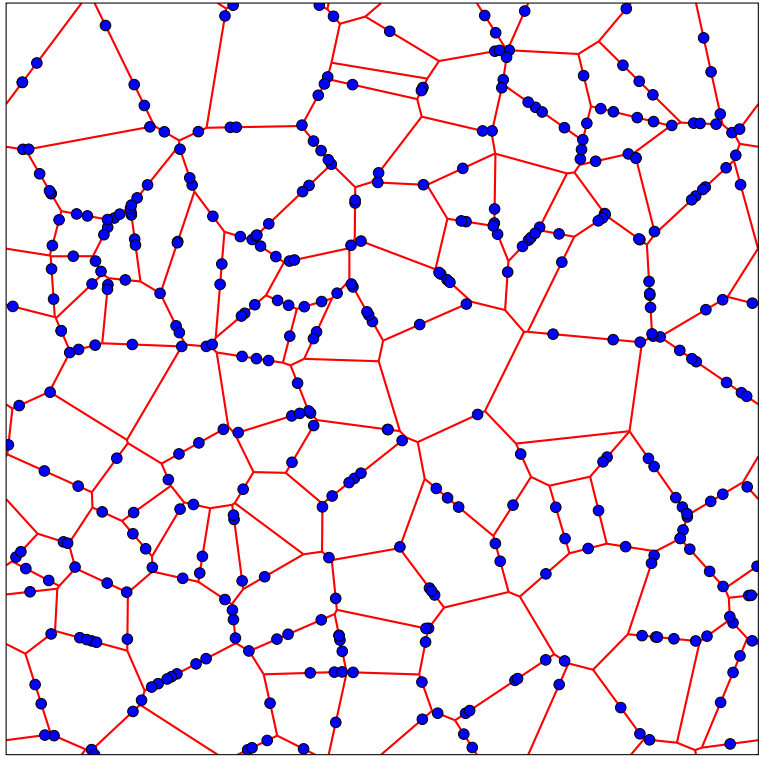}
\caption{Realization of initial device positions (blue) confined to a street system given by a Poisson--Voronoi tessellation.}
\label{Fig_Urban_Points}
\end{figure}

The second model layer consists of the point process $X^\la = \{X_i\}_{i\ge1}$ representing {\em initial positions of  devices} scattered at random in $\R^2$. 
%Here, we assume that the positions in $X^\la$ are placed at random on the street system. 
More precisely, we assume $X^\la$ to be a {\em Cox point process} with random intensity measure
$$\Lambda_{S}(A) = \la |S\cap A|$$
for every measurable $A\subset \R^2$. Here, $\la\ge 0$ is a scaling parameter that allows us to tune the expected number of devices per unit of street length, see Figure~\ref{Fig_Urban_Points} for an illustration. 

Note that, $X^\la$ can be seen an independent Poisson point processes with parameter $\la$ confined to the environment $S$. In particular for any stationary street system $S$, also the associated Cox point process of initial device positions $X^\la$ is stationary. %Recall then that the parameters $\la$ and $\gamma$ can be used to tune the model for example towards an urban setting, where initial device positions and streets are rather dense, or to a rural setting, where devices are sparsely distributed along a sparse street system. Note that the intensity of devices per unit area is simply given by $\la\g$. %\red{GAMMA}

%%%%%%%%%%%%%%%%%%%%%%%%%%%%%%%%%%%%%%%%%%%%%%%%%%%%
%Mobility
%%%%%%%%%%%%%%%%%%%%%%%%%%%%%%%%%%%%%%%%%%%%%%%%%%%%
\subsection{Device mobility via random-waypoint models}\label{mob_sec}
The third layer of the model concerns device mobility.  We present an augmented version of the classical random-waypoint model, in which a device moves at constant speed and visits a sequence of waypoints distributed in space~\cite{bettHart}. This model describes basic features of device mobility, but neglects the important constraint that in cities, devices can only move along the street system. 
For this, we assume that each initial device position $x\in S$ picks a target location $y\in S$ independently at time zero, using the probability kernel
probability kernel
$$
\k^S(x, {\rm d} y)
$$
that respects the street system $S$. In the following, we will always assume that $\k$ is {\em translation covariant}, i.e., for all $z\in \R^2$,
probability kernel
$$
\k^S(x, {\rm d} y)=\k^{S+z}(x+z, {\rm d} y+z).
$$
Examples that we will use for this {\em waypoint kernel} $\k$ are
$$
\k'^{S}(x,{\rm d} y)=|B_L(x)|^{-1}|q_S^{-1}({\rm d} y\cap S)\cap B_L(x)|\quad\text{or}\quad \k''^{S}(x,{\rm d} y)=|S\cap B_L(x)|^{-1}|{\rm d} y\cap S\cap B_L(x)|, 
$$
where $q_S(y)$ denotes the closest point on $S$ to $y$ and $B_L(x)$ denotes the disc of radius $L$, centered at $x\in \R^2$. In words, in the kernel $\k'$ the target location is chosen to be the point on the street system that is closest to a uniformly chosen point within a disc with radius $L>0$ around the location of the device. In the kernel $\k''$ the target location is chosen uniformly on the street system within a disc around the position of radius $L>0$. Note that both $\k'$ and $\k''$ are probability kernels in the sense that for almost-all $S$ and all $x\in S$, $\k'^{S}(x,\R^2)=\k''^{S}(x,\R^2)=1$. The {\em support} of $\k^S(x,\d y)$ is defined via
$$\supp(\k^S(x,\d y))=\{y\in S\colon \k^S(x,B_\e(y))>0\text{ for all }\e>0\},$$
and, roughly speaking, contains all reachable points on the street system from the position $x\in S$. We say that $\k$ is of {\em bounded support}, if there exists $L'>0$ such that $\supp(\k^S(x,\d y))\subset B_{L'}(x)$ for almost-all $S$ and all $x\in S$. Note that $\k''$ is of bounded support with $L'=L$ and also $\k'$ is of bounded support with $L'=2L$ for street systems of Poisson--Voronoi type since the projection $q_S(y)$ of any point $y\in B_L(x)$ towards the street system can not be at distance larger than $|x-y|$. 

Next, we equip each device $X_i\in X^\la$ with a velocity $V_i$ that is an independent random variable drawn from a general mutual  distribution $\muv$. 
In the simplest case, we can set $\muv=\delta_v$, where $v\ge 0$ then represents a globally fixed velocity. Slightly more generally, $V_i$ can have only two possible values, representing a pedestrian speed $v_{\rm p}\ge 0$ and a driving speed $v_{\rm d}\ge 0$, appearing with a probability $p$ and $1-p$, respectively, i.e., $\muv=p\delta_{v_{\rm p }} + (1-p)\delta_{v_{\rm d }}$. In the most general case, we assume that $0<v_{\rm min}\le v_{\rm max}<\infty$, where $v_{\rm min}:=\sup\{v\colon \muv([v,\infty))=1\}$ and $v_{\rm max}:=\inf\{v\colon \muv([0,v])=1\}$.

\begin{figure}[!htpb]
\centering
 \input{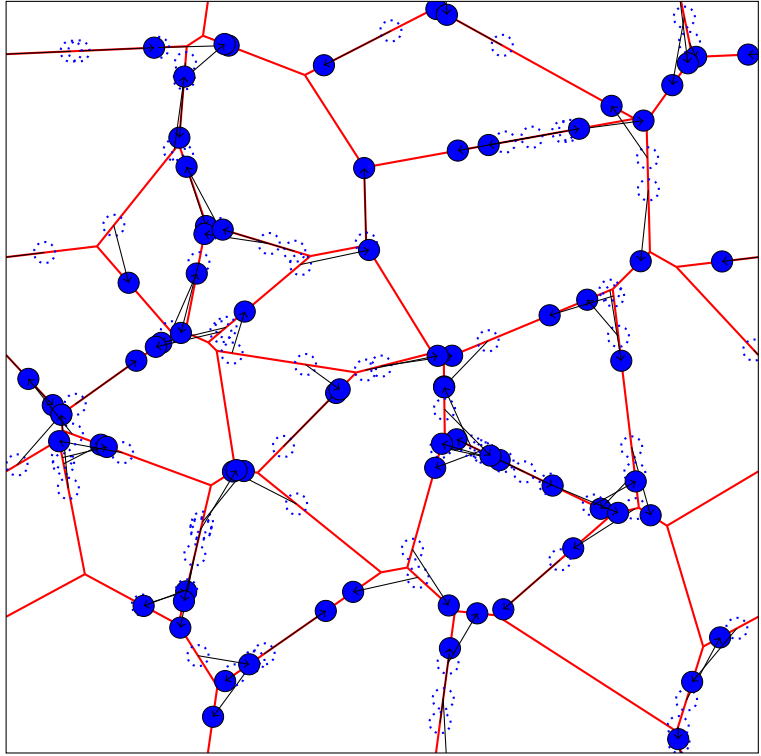}
\caption{Realization of initial device positions (dotted blue) confined to a street system given by a Poisson--Voronoi tessellation and their respective positions at a fixed positive time (blue), with arrows indicating the corresponding displacement.}
\label{Fig_Mobil_Points}
\end{figure}

Then, we assume that each device $X_i$ moves towards its individual target $Y_i$ with speed $V_i$ on the shortest path along the street system. As $S$ is asymptotically essentially
connected, $X_i$ and $Y_i$ lie in the same component and such a connecting path exists. Upon arrival at the target location, each device immediately returns to its initial position using the same path and velocity and subsequently repeats this commute until the end of the mutual time horizon $T>0$, see Figure~\ref{Fig_Mobil_Points} for an illustration. 

On the general level, we will use the notation
$$\k_t^S(x,\d y),$$
for the {\em transition kernel} that describes the probability to see the devices with initial position at $x\in\R^2$ at position $\d y\in\R^2$ at time $t$, depending on the street-system realization $S$ as well as the choice of the target location under the waypoint kernel. 
In order to simplify the notation, similarly to the street system, we are going to suppress the dependence of the theorems on the movement kernel wherever it is possible.%Let us introduce and verify some properties of our transition kernels, partially first presented in~\cite[Report II]{CRE4}. 

\subsection{Connectivity in mobile D2D networks via space-time vicinity rules}\label{con_sec}
We have now defined the individual movement patterns of the devices. Let us abbreviate by $\tra_i=(\tra_{i,t})_{0\le t\le T}$ the whole trajectory of device $X_i$ up to time $T$ based on the scheme just described, where $\tra_{i,t}$ denotes the position of device $X_i$ at time $t$, started at the initial position $\tra_{i,0}=X_i$. We now turn to the description of a connection between device $X_i$ and $X_j$. For this, we will call any intersection point of two distinct segments in $S$ or the end of a segment a {\em crossing}. Two distinct crossings are considered to be neighbors if and only if they are associated to the same segment and not separated by another crossing. Any part of a segment that lies between two neighboring crossings is called a {\em street}. Using these definitions, we let
$$
Z(X_i,X_j)=\{t\ge 0\colon |\tra_{i,t}-\tra_{j,t}|<r\text{ and } \tra_{i,t}, \tra_{j,t}\text{ are on the same street}\}
$$
denote the set of {\em contact times} between the devices, i.e., the set of times at which the devices associated with the initial positions $X_i$ and $X_j$ are on the same street and close together, where closeness is defined {via the euclidean distance} being smaller than a global parameter $r>0$. The condition that connections are only possible if devices are on the same street reflects the {\em shadowing} phenomenon, i.e., the fact that transmissions through houses is effectively blocked. Let us note that such a constraint was already introduced and investigated in~\cite{GBCE19a,GBCE19}, although in a non-dynamic situation. 

\begin{figure}[!htpb]
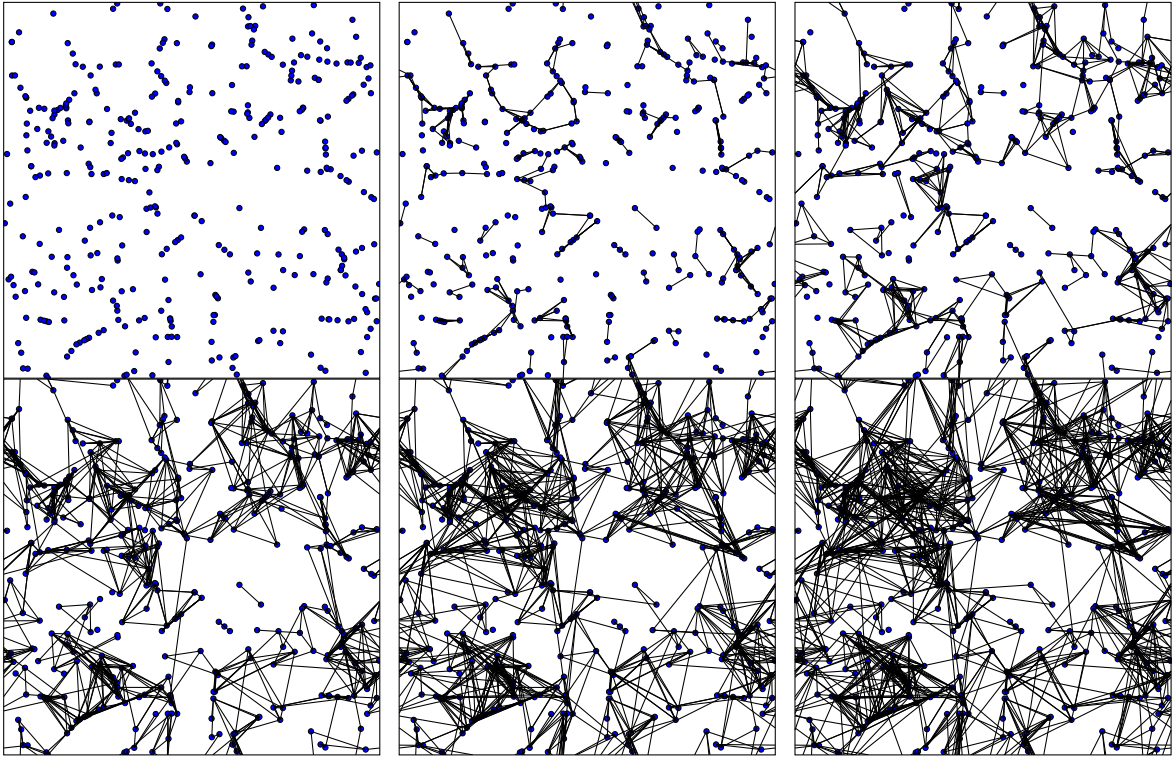

\centering
 \input{Picture/pic5_connections_T=0.tex}
 \input{Picture/pic5_connections_T=1.tex}
 \input{Picture/pic5_connections_T=2.tex}

 \input{Picture/pic5_connections_T=3.tex}
 \input{Picture/pic5_connections_T=4.tex}
 \input{Picture/pic5_connections_T=5.tex}

\caption{Realization of initial device positions (blue, up and left) confined to a street system (visually suppressed) and the appearance of edges (black lines) over time (up and left to down and right).}
\label{Fig5}
\end{figure}

Further, let $\rho=\rho_o+\rho_1$ be the sum of a fixed initialization time $\rho_o\ge 0$ and a fixed transmission time $\rho_1>0$. In other words, $\rho$ represents the minimal amount of time that two devices have to be in contact with each other in order to successfully complete a data transmission. Now, using these definitions, we say that the devices $X_i$ and $X_j$ are connected, represented by an {\em edge} in the random graph with vertices $X^\la$, if there exists $T\ge t\ge 0$ such that 
\begin{align}\label{eq_connectedness}
[t,t+\rho]\subset Z(X_i,X_j), 
\end{align}
see Figure~\ref{Fig5} for an illustration of the resulting space-time connectivity network.
%
%
%%====================================================
%%====================================================
%%====================================================
%
\section{Results}\label{sec_results}
%	- Existence of sub & super critical regimes
%	- percolation probabilities
%	- pair connection probabilities
%	- Selfavoiding paths
In the previous section we have described a multi-parameter connectivity network and the set of results presented in this article evolve around its clustering behavior. More precisely, we consider the associated random graph of connections up to time $T$
$$g_{T,\muv,\rho,r}(X^\la)=g_{\k,T,\mathcal{S},\muv,\rho,r}(X^\la),$$
in which we will suppress the dependence on the street-system as well as the mobility kernel $\kappa$, since we will mostly assume those parameters to be fixed. Recall that $g_{T,\muv,\rho,r}(X^\la)$ is the graph with vertices given by the Cox points in $X^\la$ and in which any pair of Cox points is connected if~\eqref{eq_connectedness} holds for at least one $t \in [0, T-\rho]$, i.e., their associated trajectories (described via $\k$ and $\muv$) are sufficiently close together (indicated by the parameter $r$) on the same street, for a sufficient amount of time (indicated by the parameter $\rho$). Let us highlight that, although the underlying mechanism for the edge-drawing depends on time, the graph $g_{T,\muv,\rho,r}(X^\la)$ is static, see also below for more comments. 

In view of the propagation of malware, or more generally any data, through the system, we start our analysis by considering percolation properties of $g_{T,\muv,\rho,r}(X^\la)$. That is, we say that the graph {\em percolates} if it contains an unbounded connected component. Here, percolation can be seen as a certain best-case scenario for the dissemination of data in the system. Indeed, our notion of percolation does not reflect that connections are in fact established at specific points in time, and that the actual path of a message through the dynamical network depends on when the message can hop from one device to another. In that sense, from the perspective of a message in the system, percolation assumes an optimal space-time placement of the connection events. In other words, in the absence of percolation, even an optimal space-time placement of the connection events will not lead to the possibility that data travels very far. 

Let us mention that a percolation analysis can additionally be motivated by assuming that the devices travel along their trajectory repeatedly (e.g., every day). Then, in case the handover of a message over multiple hops is impossible in one time window $[0,T]$, in the second iteration, when the possible transmissions in the first iteration are already performed, another set of transmissions will happen. 
Repeating this mechanism ultimately leads to a complete coverage of the connected component and then percolation indicates whether or not this component is large enough to support long-distance transmissions. 

%To explain these motivations, let us make an example. Imagine three devices $X_1,X_2$ and $X_3$, where we see only a connection event at time $t_1\in [0,T]$ between $X_1$ and $X_2$ and another connection event at time $t_2\in [0,T]$ between $X_2$ and $X_3$. In our percolation picture, all three devices belong to the same connected component independent of the ordering of the times. However, if we entertain the idea that device $X_1$ sends a message to the devices in its component, in particular, to device $X_3$, then, in one time window this would require $t_1<t_2$, where then $X_2$ serves as a relay. Hence, in our first motivation, only in a connected component devices can exchange messages, however, we suppress the additional requirement that there also needs to be a sequence of connection events in the right order, in order to send messages in one time window. In our second motivation, which becomes relevant for example when dealing with malware, coming back to our example, even though $t_1<t_2$, in the first time window the malware can jump from $X_1$ to $X_2$ and then, in the second time window finally arrives at $X_3$. 

Consider for $a>0$ the following rescaled distributions, 
\begin{itemize}
\item $\k^{S}[a](x, \d y):=\kappa^{S/a}(x/a,\d y/a)$,
\item $\mathcal{S}^a(\d S):=\mathcal{S}(\d S/a)$, and 
\item $\muv^{a}(\d v):=\muv(\d v/a)$,
\end{itemize}
%let us denote by $\k^{S}[a](x, \d y)=\kappa^{S/a}(x/a,\d y/a)$ the spatially re-scaled kernel and $\mathcal{S}^a$ the re-scaled version of the street system with $\mathcal{S}^a(\d S)=\mathcal{S}(\d S/a)$.
%Similarly, we define $\muv^a$ to be the re-scaled velocity via
%$$\int\muv^{a}(\d v)f(v)=\int\muv(\d v)f(av),$$
%for all bounded and measurable functions $f: [0,\infty)\to\R$. 
%In words, $\muv^{a}$ describes a model in which devices move $a$-times as fast as in the original model.
then, we have the following scaling relations with respect to percolation for which we will present the short proof in Section~\ref{proof_prop_scale}. 
\begin{proposition}\label{prop_scale}
For all $a>0$ we have that 
\begin{enumerate}
\item[(1)] $g_{\k[a],T,\mathcal{S}^{1/a},\muv^a,\rho,ar}(X^{\la/a})$ percolates if and only if $g_{\k,T,\mathcal{S},\muv,\rho,r}(X^\la)$ percolates and
\item[(2)] $g_{\k,T/a,\mathcal{S},\muv^a,\rho/a,r}(X^{\la})$ percolates if and only if $g_{\k,T,\mathcal{S},\muv,\rho,r}(X^\la)$ percolates.
\end{enumerate}
\end{proposition}
Let us now fix $\k$ and $\mathcal S$ and define the following critical parameters for percolation,
\begin{align*}
\Tc=\Tc(\la,\muv,\rho,r)&:=\inf\{\Tc>0\colon \P(g_{T,\muv,\rho,r}(X^\la)\text{ percolates})>0\}\\
\lac=\lac(T,\muv,\rho,r)&:=\inf\{\la\ge0\colon \P(g_{T,\muv,\rho,r}(X^\la)\text{ percolates})>0\}\\
\rhoc=\rhoc(T,\la,\muv,r)&:=\sup\{\rho\ge 0\colon \P(g_{T,\muv,\rho,r}(X^\la)\text{ percolates})>0\}\\
\rc=\rc(T,\la,\muv,\rho)&:=\inf\{r\ge0\colon \P(g_{T,\muv,\rho,r}(X^\la)\text{ percolates})>0\}
\end{align*}
%\begin{align*}
%\Tc&:=\inf\{\Tc>0\colon \P(g_{\k,T,\g,\muv,\rho,r}(X^\la)\text{ percolates})>0\}\\
%\lac&:=\inf\{\la>0\colon \P(g_{\k,T,\g,\muv,\rho,r}(X^\la)\text{ percolates})>0\}\\
%\rhoc&:=\sup\{\rho>0\colon \P(g_{\k,T,\g,\muv,\rho,r}(X^\la)\text{ percolates})>0\}\\
%\rc&:=\inf\{r>0\colon \P(g_{\k,T,\g,\muv,\rho,r}(X^\la)\text{ percolates})>0\}
%\end{align*}
and note that these critical thresholds describe subcritical regimes in the sense that for $T<\Tc$, $\la<\lac$, $\rho>\rhoc$ or $r<\rc$, the system features percolation with probability zero. Further note that percolation is indeed monotone with respect to these parameters since longer time horizons, more devices, shorter connection-time thresholds and larger vicinity thresholds can only increase connectivity. Hence, we will present conditions under which $0<\Tc,\lac,\rhoc,\rc<\infty$. 

However, considering the dependence of percolation on the velocity $\muv$, the situation is substantially more complicated. Indeed, keeping all other parameters fixed and assuming for example a small density of devices, then, slow velocities may not allow any message to travel very far, and thus we expect to see absence of percolation. On the other hand, if velocities are large, then large gaps between devices can be bridged, but at the same time, the probability to be in close proximity for a sufficiently long time becomes small, again leading to absence of percolation. In an intermediate speed regime in which also the other parameters allow for percolation, we then expect to see percolation with positive probability. 

Let us recall that we generally assume the support of $\muv$ to be contained in the minimal interval $[v_{\rm min},v_{\rm max}]$ with $0<v_{\rm min}\le v_{\rm max}<\infty$. 
We hence define two critical speed parameters
\begin{align*}
\vc=\vc(T,\la,\rho,r)&:=\inf\{av_{\rm max}\ge 0\colon \P(g_{T,\muv^a,\rho,r}(X^\la)\text{ percolates})>0\}\\
\vcc=\vcc(T,\la,\rho,r)&:=\sup\{av_{\rm min}\ge0\colon \P(g_{T,\muv^a,\rho,r}(X^\la)\text{ percolates})>0\}
\end{align*}
and present conditions under which $0<\vc<\vcc<\infty$.  

\medskip
In the following Section~\ref{sec-noperc} we present our results on absence of percolation and in Section~\ref{sec-perc} we present our results on percolation regimes.  
%Let us start by investigating regimes where percolation is impossible.  
\subsection{Regimes for absence of percolation}\label{sec-noperc}
Our first result establishes that indeed, insufficiently many devices and a too slow speed lead to absence of percolation. 
\begin{theorem}\label{thm_lowintense}
For any $\kappa$, $\mathcal S$ and all $r,\rho \ge 0$ the following holds:
\begin{enumerate}
\item[(1)] For all $\la\ge 0$ and $\muv$ we have that $\Tc(\la,\muv,\rho,r)>2\rho$.
\item[(2)] For all $T\ge 0$ and $\muv$ we have that $\lac(T,\muv,\rho,r)>0$.
\item[(3)] For all $T,\la\ge 0$ we have that $\vc(T,\la,\rho,r)>0$. 
\end{enumerate}
\end{theorem}
% \begin{proof}
The proof of (2) rests on the idea that we can bound the process by the easier-controllable Cox--Gilbert graph. %$g_{2(T-\rho)v_{\rm max}+r}(X^\la)$ and prove (2).
For (1) and (3) this comparison fails, and we need to couple our space-continuous model to site-percolation in $\Z^d$ and employ stabilization arguments. We present the detailed proof in Section~\ref{Sec_Proof_Sub_Rho}.

% Each device can travel at most to a distance $Tv_{\rm max}$ from its starting position. Hence, any connected component of $g_{\k,T,\g,\muv,\rho,r}(X^\la)$ is contained in a connected component of the Gilbert graph $g_{2(T-\rho)v_{\rm max}+r}(X^\la)$. In the Gilbert graph, any pair of points is connected if and only if their mutual distance is smaller than ${2(T-\rho)v_{\rm max}+r}$. As the devices need time $\rho$ to form a connection, the last connection attempt can start at time $T-\rho$ if the devices are within distance $r$ of each other.
%  But $g_{2(T-\rho)v_{\rm max}+r}(X^\la)$ percolates with probability zero for all sufficiently small  $\la$, see, for example~\cite{wias8}. The same argument only applies for $T$ and $v_{\rm max}$ for $r = 0$ (and $\rho = 0$ as a consequence, see Theorem~\ref{Thm_Rho_0}). For $r > 0$, it is possible that the Gilbert graph $g_{r}(X^\la)$ percolates. Therefore, other arguments that rely on the street system are required. In Section~\ref{Sec_Proof_Sub_Rho} we will complete the proof, that either for $v_{\rm max} \downarrow 0$ or $T \downarrow 2 \rho$ no percolation occurs as the probability for a connection to form goes to zero.
% \end{proof}
Let us mention that the comparison to Gilbert graphs also provides rough upper bounds for $\lac(T,\muv,\rho,r)$ and $\vc(T,\la,\rho,r)$ with the help of the numerical analysis provided in \cite{wias9,wias2}. 

\medskip
Let us mention that item (3) in Theorem~\ref{thm_lowintense} does not hold for $T=\infty$. In this case we have $\vc(\infty,\la,\rho,r)=0$ once $\P(g_{\infty,\muv,\rho,r}(X^\la)\text{ percolates})>0$ for some parameter configuration. Indeed, using Part (2) in Proposition~\ref{prop_scale}, we also have percolation for all $0<a<1$, i.e.,  $0<\P(g_{\infty,\muv^{a},\rho/a,r}(X^\la)\text{ percolates})\le \P(g_{\infty,\muv^{a},\rho,r}(X^\la)\text{ percolates})$, where we also used monotonicity in $\rho$.

\medskip
Next, note that $\rc$ can very well be zero, since devices with non-parallel movement on a street share the same position exactly once, leading to a connection even for $r=0$, if $\rho=0$. However, devices that move together in the same direction can not exchange messages if $r=0$. For the consequences of this observation regarding the critical connectivity threshold, see Theorem~\ref{Thm_Rho_1} below. Let us note that, in principle, absence of percolation for $r=0$ could still imply $\rc(T,\la,\muv,\rho)=0$. This would be the case if $g_{T,\muv,\rho,r}(X^\la)$ percolates for all positive $r$, which indicates that the parameters $(T,\la,\muv,\rho,0)$ are critical for percolation, rendering a mathematical analysis extremely hard.

\medskip
In order to see absence of percolation for large speeds $v_{\rm min}$, first note that any street $s$ such that $2|s|/v_{\rm min}<\rho$ can not be used to exchange messages. However, every street can be used for traveling. More precisely, in order to establish a connection, any device has to remain at least a time of $\rho$ on some street and during that time travels at least a distance of $\rho v_{\rm min}$. As devices reverse their direction when they reach either their destination or their initial position, the minimal length of a street such that communication on it is possible is given by $\rho v_{\rm min}/2$. Although staying on a street of length smaller than $\rho v_{\rm min}/2$ is possible if both the initial position and the destination are on that street, those devices can only connect to the same type of devices that have both their initial position and destination on the same street. In that way, only local connective components are formed, we will ignore this possibility in the following considerations.

We have the following result on absence of percolation for sufficiently high velocities. 
\begin{theorem}\label{thm_largespeeds}
We consider street systems of Poisson--Voronoi or Poisson--Delaunay type.
%, which are in particular exponentially stabilizing and asymptotically essentially connected.  
Then, for any $\k$ with bounded support and all $T, \la,r,\rho\ge 0$ we have that $\vcc(T,\la,\rho,r)<\infty$. 
\end{theorem}
Let us note that we in fact prove that $\vcc(\infty,\la,0,r)<\infty$, which implies a uniform bound in $T$ and $\rho$ by the corresponding monotonicites.

The proof rests on the following ideas. As $v_{\rm min}$ tends to infinity, any device $X_i\in S$ eventually reaches its target location $Y_i\in S$ before time $T$.
%, where $|X_i-Y_i|\le \diam(\supp(\k^S(X_i,\d y)))$.
In particular, it is very unlikely to see a path from $X_i$ to any of its possible target locations that reaches very far away from $X_i$ in $\R^2$, and these shortest paths can be determined independently of the device process. Additionally, only devices that have a sufficiently long street on a shortest path to at least one of their possible target locations can contribute to the connectivity, and again this can be determined independently of the device process. However, as $v$ tends to infinity, existence of these long streets becomes also unlikely, leading to a thinning of the devices and eventually to a subcritical regime. 
Formally, we dominate our percolation model by a subcritical {\em geostatistical Cox--Boolean model}, see also~\cite{AT17,SRPD04}, where the (random) environment determines both, the random radii associated to the points, and a sufficiently strong thinning of the Cox points. 
We present the complete proof in Section~\ref{Sec_Proof_Sub}.

\medskip
We still have to deal with the connection-time parameter $\rho$. For this, note that we have a trivial bound for the connection time, given by $\rhoc\le T/2$ since for $\rho> T/2$, %there can never be an edge between initial device positions on different streets.
in order to establish a connection on at least two different streets, the device has to stay on one street for at least $T/2$ time, leaving not enough time to form connections on any different street. Therefore, any device can only have connections to other devices that also spend the majority of their time on the same street.  As a consequence, each street forms its local connection cluster and communication between different clusters is impossible, which implies absence of percolation due to the stabilization assumption.
However, this trivial bound can be substantially improved, see below. 

Let us first make a statement on the absence of percolation for any positive $\rho$ if $r= 0$. As almost surely devices do not share the same position while moving in the same direction, communication is only possible between devices that move in opposite direction and $\rho=0$ becomes a necessary condition.
Therefore
\begin{align*}
    \rhoc(T,\la,\muv,0)\in \{0,-\infty\},
\end{align*}
since even with $\rho = 0$, percolation might not be possible if the device density $\la$ is too small.
Interestingly, $\rho = 0$ is the only case where the model becomes monotone in $v$, as the distance a device travels during a successful communication $\rho v = 0$ does not increase in $v$.

% \blue{Let us first make a statement on the absence of percolation for any positive $\rho$. 
% \begin{remark}\label{Thm_Rho_0}
% For any $\k$ and all $T,\g,\la\ge 0$ we have that $\rhoc(\k,T,\g,\la,\muv,0)=0$.
% \end{remark}
% \begin{proof}
% As indicated above, connections between devices that move together in the same direction on a street, can not exchange messages if $r=0$. Additionally, if $r=0$, devices in non-parallel movement can never have a contact time that includes intervals of positive length, and hence $g_{\k,T,\g,\muv,\rho,0}(X^\la)$ does not percolate with positive probability for any $\rho>0$. 
% \end{proof} }

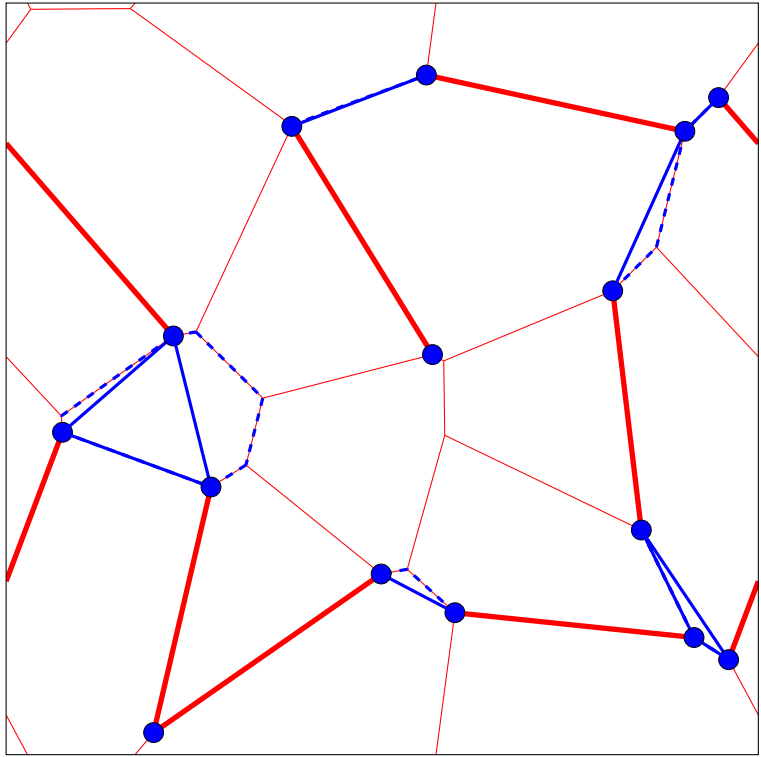
\begin{figure}[h]
\centering
\begin{tikzpicture}[scale=1.25] 
 \begin{scope} 
\clip(-4,-4) rectangle (4,4);
\draw[red, line width = 2pt](0.5346409023673402,0.258314289366571)--(-0.9613066390742063,2.687401972219394);
\draw[red, thin](0.5346409023673402,0.258314289366571)--(-1.2700999778371633,-0.20535476335735314);
\draw[red, thin](0.5346409023673402,0.258314289366571)--(0.6545503199713438,0.19147992910574074);
\draw[red, thin](-0.9613066390742063,2.687401972219394)--(-1.981262656708152,0.5000831257795172);
\draw[red, thin](-0.9613066390742063,2.687401972219394)--(-0.9631031311463767,2.6999087950218934);
\draw[red, thin](-1.2700999778371633,-0.20535476335735314)--(-1.981262656708152,0.5000831257795172);
\draw[red, thin](-1.2700999778371633,-0.20535476335735314)--(-1.449304249986116,-0.9165803028143067);
\draw[red, thin](-1.981262656708152,0.5000831257795172)--(-2.2212212080044775,0.45640537163588046);
\draw[red, thin](-3.4132148866090684,-0.3945847551872484)--(-3.400567306034886,-0.5685512266627061);
\draw[red, thin](-3.4132148866090684,-0.3945847551872484)--(-2.2212212080044775,0.45640537163588046);
\draw[red, thin] (-3.4132148866090684,-0.3945847551872484)--(-4,0.2362372426536774);
\draw[red, thin] (4,0.2362372426536774)--(2.9162579606549235,1.4013116503072265);
\draw[red, thin](-3.400567306034886,-0.5685512266627061)--(-1.821173486341278,-1.1510917968917462);
\draw[red, line width = 2pt] (-3.400567306034886,-0.5685512266627061)--(-4,-2.153568907805884);
\draw[red, line width = 2pt] (4,-2.153568907805884)--(3.6840554171198785,-2.9889883888719986);
\draw[red, line width = 2pt] (-2.2212212080044775,0.45640537163588046)--(-4,2.50638745327403);
\draw[red, line width = 2pt] (4,2.50638745327403)--(3.5772305633663386,2.9936148553901685);
\draw[red, thin](-1.449304249986116,-0.9165803028143067)--(-1.821173486341278,-1.1510917968917462);
\draw[red, thin](-1.449304249986116,-0.9165803028143067)--(-0.011401148090081581,-2.0772586447779116);
\draw[red, line width = 2pt](-1.821173486341278,-1.1510917968917462)--(-2.4313590321486274,-3.7649909492486575);
\draw[red, line width = 2pt](2.4515058294811296,0.9372423847721281)--(2.755682905768952,-1.608975122414948);
\draw[red, thin](2.4515058294811296,0.9372423847721281)--(0.6545503199713438,0.19147992910574074);
\draw[red, thin](2.4515058294811296,0.9372423847721281)--(2.9162579606549235,1.4013116503072265);
\draw[red, thin](2.755682905768952,-1.608975122414948)--(0.6650867647420391,-0.600419400157248);
\draw[red, thin](2.755682905768952,-1.608975122414948)--(3.316340046926652,-2.7532264881417214);
\draw[red, thin](0.6545503199713438,0.19147992910574074)--(0.6650867647420391,-0.600419400157248);
\draw[red, thin](0.6650867647420391,-0.600419400157248)--(0.2661661979260624,-2.026423027609807);
\draw[red, thin](2.9162579606549235,1.4013116503072265)--(3.218742062269908,2.6351233652626274);
\draw[red, thin](3.6840554171198785,-2.9889883888719986)--(3.316340046926652,-2.7532264881417214);
\draw[red, thin] (3.6840554171198785,-2.9889883888719986)--(4,-3.577215445794999);
\draw[red, thin] (-4,-3.577215445794999)--(-3.772916787736368,-4);
\draw[red, thin] (-3.772916787736368,4)--(-3.7372385963767405,3.93357418155901);
\draw[red, line width = 2pt](3.316340046926652,-2.7532264881417214)--(0.7729427770792394,-2.489081205800572);
\draw[red, thin](-0.9631031311463767,2.6999087950218934)--(0.4697231481860521,3.233684719679209);
\draw[red, thin](-0.9631031311463767,2.6999087950218934)--(-2.6780737982759124,3.9411650797759084);
\draw[red, line width = 2pt](0.4697231481860521,3.233684719679209)--(3.218742062269908,2.6351233652626274);
\draw[red, thin] (0.7729427770792394,-2.489081205800572)--(0.5717600171728314,-4);
\draw[red, thin] (0.5717600171728314,4)--(0.4697231481860521,3.233684719679209);
\draw[red, thin](3.218742062269908,2.6351233652626274)--(3.5772305633663386,2.9936148553901685);
\draw[red, thin](-2.6780737982759124,3.9411650797759084)--(-3.7372385963767405,3.93357418155901);
\draw[red, thin] (-2.4313590321486274,-3.7649909492486575)--(-2.6286753248785537,-4);
\draw[red, thin] (-2.6286753248785537,4)--(-2.6780737982759124,3.9411650797759084);
\draw[red, thin] (-3.7372385963767405,3.93357418155901)--(-4,3.5732912921702065);
\draw[red, thin] (4,3.5732912921702065)--(3.5772305633663386,2.9936148553901685);
\draw[red, thin](0.7729427770792394,-2.489081205800572)--(0.2661661979260624,-2.026423027609807);
\draw[red, thin](0.2661661979260624,-2.026423027609807)--(-0.011401148090081581,-2.0772586447779116);
\draw[red, line width = 2pt](-2.4313590321486274,-3.7649909492486575)--(-0.011401148090081581,-2.0772586447779116);
\draw[fill = blue ,ultra thin] (0.5346409023673402,0.258314289366571) circle(3pt) ;
\draw[fill = blue ,ultra thin] (-0.9613066390742063,2.687401972219394) circle(3pt) ;
\draw[blue, very thick](-0.9613066390742063,2.687401972219394)--(0.4697231481860521,3.233684719679209);
\draw[blue, dashed,very thick](-0.9613066390742063,2.687401972219394)--(-0.9631031311463767,2.6999087950218934);
\draw[blue, dashed,very thick](-0.9631031311463767,2.6999087950218934)--(0.4697231481860521,3.233684719679209);
\draw[fill = blue ,ultra thin] (-3.400567306034886,-0.5685512266627061) circle(3pt) ;
\draw[blue, very thick](-3.400567306034886,-0.5685512266627061)--(-1.821173486341278,-1.1510917968917462);
\draw[blue, dashed,very thick](-3.400567306034886,-0.5685512266627061)--(-1.821173486341278,-1.1510917968917462);
\draw[blue, very thick](-3.400567306034886,-0.5685512266627061)--(-2.2212212080044775,0.45640537163588046);
\draw[blue, dashed,very thick](-3.400567306034886,-0.5685512266627061)--(-3.4132148866090684,-0.3945847551872484);
\draw[blue, dashed,very thick](-3.4132148866090684,-0.3945847551872484)--(-2.2212212080044775,0.45640537163588046);
\draw[fill = blue ,ultra thin] (2.4515058294811296,0.9372423847721281) circle(3pt) ;
\draw[blue, very thick](2.4515058294811296,0.9372423847721281)--(3.218742062269908,2.6351233652626274);
\draw[blue, dashed,very thick](2.4515058294811296,0.9372423847721281)--(2.9162579606549235,1.4013116503072265);
\draw[blue, dashed,very thick](2.9162579606549235,1.4013116503072265)--(3.218742062269908,2.6351233652626274);
\draw[fill = blue ,ultra thin] (2.755682905768952,-1.608975122414948) circle(3pt) ;
\draw[blue, very thick](2.755682905768952,-1.608975122414948)--(3.316340046926652,-2.7532264881417214);
\draw[blue, dashed,very thick](2.755682905768952,-1.608975122414948)--(3.316340046926652,-2.7532264881417214);
\draw[blue, very thick](2.755682905768952,-1.608975122414948)--(3.6840554171198785,-2.9889883888719986);
\draw[blue, dashed,very thick](2.755682905768952,-1.608975122414948)--(3.316340046926652,-2.7532264881417214);
\draw[blue, dashed,very thick](3.316340046926652,-2.7532264881417214)--(3.6840554171198785,-2.9889883888719986);
\draw[fill = blue ,ultra thin] (3.5772305633663386,2.9936148553901685) circle(3pt) ;
\draw[blue, very thick](3.5772305633663386,2.9936148553901685)--(3.218742062269908,2.6351233652626274);
\draw[blue, dashed,very thick](3.5772305633663386,2.9936148553901685)--(3.218742062269908,2.6351233652626274);
\draw[fill = blue ,ultra thin] (3.6840554171198785,-2.9889883888719986) circle(3pt) ;
\draw[blue, very thick](3.6840554171198785,-2.9889883888719986)--(3.316340046926652,-2.7532264881417214);
\draw[blue, dashed,very thick](3.6840554171198785,-2.9889883888719986)--(3.316340046926652,-2.7532264881417214);
\draw[fill = blue ,ultra thin] (-2.4313590321486274,-3.7649909492486575) circle(3pt) ;
\draw[fill = blue ,ultra thin] (3.316340046926652,-2.7532264881417214) circle(3pt) ;
\draw[fill = blue ,ultra thin] (0.7729427770792394,-2.489081205800572) circle(3pt) ;
\draw[blue, very thick](0.7729427770792394,-2.489081205800572)--(-0.011401148090081581,-2.0772586447779116);
\draw[blue, dashed,very thick](0.7729427770792394,-2.489081205800572)--(0.2661661979260624,-2.026423027609807);
\draw[blue, dashed,very thick](0.2661661979260624,-2.026423027609807)--(-0.011401148090081581,-2.0772586447779116);
\draw[fill = blue ,ultra thin] (-0.011401148090081581,-2.0772586447779116) circle(3pt) ;
\draw[fill = blue ,ultra thin] (-2.2212212080044775,0.45640537163588046) circle(3pt) ;
\draw[blue, very thick](-2.2212212080044775,0.45640537163588046)--(-1.821173486341278,-1.1510917968917462);
\draw[blue, dashed,very thick](-2.2212212080044775,0.45640537163588046)--(-1.981262656708152,0.5000831257795172);
\draw[blue, dashed,very thick](-1.981262656708152,0.5000831257795172)--(-1.2700999778371633,-0.20535476335735314);
\draw[blue, dashed,very thick](-1.2700999778371633,-0.20535476335735314)--(-1.449304249986116,-0.9165803028143067);
\draw[blue, dashed,very thick](-1.449304249986116,-0.9165803028143067)--(-1.821173486341278,-1.1510917968917462);
\draw[fill = blue ,ultra thin] (-1.821173486341278,-1.1510917968917462) circle(3pt) ;
\draw[fill = blue ,ultra thin] (0.4697231481860521,3.233684719679209) circle(3pt) ;
\draw[fill = blue ,ultra thin] (3.218742062269908,2.6351233652626274) circle(3pt) ;
 \end{scope} 
 \draw (-4,-4) rectangle (4,4); 
 \end{tikzpicture}
\caption{Realization of a street system given by a Poisson--Voronoi tessellation (red). Streets of length greater than $\rho v/2$ are drawn thick. Their endpoints (blue) are connected if their graph distance (indicated in dashed blue lines) is less than $v(T-2\rho)$. The blue vertices, combined with the thick red and thick blue edges form the percolation graph  $S^{v\rho/2,v(T-2\rho)}$.}
\label{Fig_long_edge_percolation}
\end{figure}

In order to provide potentially better upper bounds for the critical connection time $\rhoc$, we consider the following percolation model based on the street system. 
For all $S$ first denote by $S^{a}$ the thinned street system consisting of all streets in $S$ with length not smaller than $a\ge 0$. Second, we connect any pair of vertices that are endpoints of any of the long streets (with length not shorter than $a$) if and only if the length of the shortest path along the streets between them does not exceed $b\ge 0$. The resulting graph $S^{a,b}$ now consists of long streets and special edges between some of the endpoints of the long streets, see Figure~\ref{Fig_long_edge_percolation} for an illustration. 
We say that $S^{a,b}$ percolates, if $S^{a,b}$ contains an infinitely-long component of connected streets and edges. Note that $\rho\mapsto \one\{S^{v_{\rm min}\rho/2,v_{\rm max}(T-2\rho)}\text{ percolates}\}$ is decreasing since we see less surviving streets and also less connections between endpoints of long streets. Here the graph is defined only for $0 \le \rho \le T/2$. For $\rho > T/2$, the allowed distance $v(T-2\rho)$ becomes negative. Hence, we can define a unique critical threshold via 
\begin{align*}
\rhoc'(T,\muv)&:=\sup\{0\le\rho\le T/2\colon \P(S^{v_{\rm min}\rho/2,v_{\rm max}(T-2\rho)}\text{ percolates})>0\}.
\end{align*}
%Further let $\rhoc'(T,\g,\muv):=\sup\{\rhoc'(T,\g,v)\colon v_{\min}\le v\le v_{\rm max}\}$ denote the largest possible critical threshold in order to see percolation for all possible speeds. 
%In principle we could have that $\rhoc'(T,\g,\muv)>T/2$ since $\P(S^{v_{\rm min}T/2,0}\text{ percolates})$ may be positive for sufficiently small $v_{\rm min}$, see Theorem~\ref{Thm_Rho_2} below for a related result. However, in this case, the reference model $S^{a,b}$ does not provide better bounds. 
Let us mention that $v\mapsto \one\{S^{v\rho/2,v(T-2\rho)}\text{ percolates}\}$ is not monotone since, for example, large velocities lead to more streets being eliminated, while also leading to the possibility to travel larger distances between remaining long streets.
Even for $T$ being much larger than $2\rho$, it is this interplay between the possibility to connect on sufficiently long streets and the ability to bridge between long streets (using all streets), that makes this percolation picture interesting. 

We have the following theorem.
\begin{theorem}\label{Thm_Rho_1}
For any $\k$, $\mathcal S$, $\muv$ and all $\la, T,r\ge 0$ we have that $\rhoc(T,\la,\muv,r)\le \min\{T/2,\rhoc'(T,\muv)\}$.
\end{theorem}
Let us note that we in fact prove that $\rhoc(T,\infty,\muv,\infty)\le \min\{T/2,\rhoc'(T,\muv)\}$, which provides a uniform upper bound in $\la$ and $r$ by monotonicity. 
Let us note that, for any fixed $\lambda$, it is possible to even prove that $\rhoc(T,\la,\delta_v,r)<T/2$ with the help of a coupling towards a discrete site-percolation model similar to the one used in the proof of Theorem~\ref{thm_lowintense}, see Section~\ref{Sec_Proof_Sub_Rho}. 
However, this approach leads to $\lim_{\lambda \uparrow \infty}\rhoc(T,\la,\delta_v,r)=T/2$.
%red

On the other hand, for any $\k$, $\la$ and $r$ we have $\rhoc(T,\la,\delta_v,r)\le \rhoc'(T,\delta_v)$ due to the construction of the auxiliary model $S^{a,b}$. Thus the analysis of subcritical regimes for the auxiliary model is the main idea of the proof, which we present in Section~\ref{Sec_Proof_Sub_7}.
 
\medskip
In order to complement Theorem~\ref{Thm_Rho_1}, as a final result for this section on absence of percolation, let us isolate a regime of fast velocities such that $\rhoc'(T,\delta_v)<T/2$. Recall that this is in general a nontrivial task since $v\mapsto \one\{S^{v\rho/2,v(T-2\rho)}\text{ percolates}\}$ is not monotone. Let us denote by 
\begin{align}\label{def_abc}
a_{\rm c}(b)&:=\sup\{a>0\colon \P(S^{a,b}\text{ percolates})>0\}
\end{align}
the largest street-thinning parameter such that we see percolation with positive probability in the presence of additional edges between endpoints of long streets. 
We further write $a_{\rm c}:=\lim_{b\downarrow0}a_{\rm c}(b)$, where we note that the limit is in fact an infimum since $b\mapsto a_{\rm c}(b)$ is decreasing. Then we have the following result. 
% \blue{ OLD VERSION
% \begin{theorem}\label{Thm_Rho_2}
% For any street system $\mathcal{S}$ we have $a_{\rm c}<\infty$. Further, for any $\k$, $T,\mathcal{S},\la,r\ge 0$ and $v>4 a_{\rm c}/T$  it holds that $\rhoc'(T,\delta_v)< T/2$.
% \end{theorem}
% In particular, the conditions presented in Theorem~\ref{Thm_Rho_2} imply that $\rhoc(T,\la,\delta_v,r)\le \rhoc'(T,\delta_v)$ and analyzing absence of percolation for the auxiliary model  $S^{a,b}$ indeed provides better bounds.
% \begin{proof}
% The proof for the finiteness of $a_{\rm c}$ rests on stabilization arguments and will be presented in Section~\ref{Sec_Proof_Sub_2}. 
% %We can assume that $v_{\rm min}=v_{\rm max}$ since $\rhoc'(T,\g,\delta_{v_{\rm min}})\ge\rhoc'(T,\g,\muv)$ where $\supp(\muv)=[v_{\rm min},v_{\rm max}]$. 
% %\begin{align*}
% %\P(S^{v_{\rm min}\rho,v_{\rm max}(T-2\rho)^+}\text{ percolates})\ge\P(S^{v_{\rm min}\rho,v_{\rm min}(T-2\rho)^+}\text{ percolates}). 
% %\end{align*}
% Under the assumption that $a_{\rm c}<\infty$, note that, since $v>4 a_{\rm c}/T$, there exist $\eps,b>0$ such that $v(T/2-\eps)/2>a_{\rm c}(b)$. In particular, we have that $\P(S^{v(T/2-\eps)/2,b'}\text{ percolates})=0$ for all $b'<b$. 
% Hence, for $\rho>T/2-\eps$ and $\rho>T/2-b/(2v)$ we have also that $0=\P(S^{v\rho/2,v(T-2\rho)}\text{ percolates})$. 
% Thus, $\rhoc'(T,\delta_v)< T/2$, which finishes this part of the proof. 
% \end{proof}
% }

\begin{theorem}\label{Thm_Rho_2}
Let $\mathcal{S}$ be a street system, then $a_{\rm c}<\infty$. Further, for any $\kappa$ and $T\ge 0$ as well as $v>4 a_{\rm c}/T$  it holds that $\rhoc'(T,\delta_v)< T/2$.
\end{theorem}
The proof for the finiteness of $a_{\rm c}$ rests on stabilization arguments. 
%We can assume that $v_{\rm min}=v_{\rm max}$ since $\rhoc'(T,\g,\delta_{v_{\rm min}})\ge\rhoc'(T,\g,\muv)$ where $\supp(\muv)=[v_{\rm min},v_{\rm max}]$. 
%\begin{align*}
%\P(S^{v_{\rm min}\rho,v_{\rm max}(T-2\rho)^+}\text{ percolates})\ge\P(S^{v_{\rm min}\rho,v_{\rm min}(T-2\rho)^+}\text{ percolates}). 
%\end{align*}
Then, we can use the monotonicity of $S^{a,b}$ to prove that, if the velocity of the devices is high enough, constructing a subcritical regime for a $\rho<T/2$ is possible. We present the details in Section~\ref{Sec_Proof_Sub_2}.

\subsection{Regimes for percolation}\label{sec-perc}
In order to see percolation we have a variety of parameters for possible adjustments. In view of Theorem~\ref{thm_lowintense}, note that for too small device intensity or speed, we can never see percolation. Hence, we assume the speed to be positive, i.e., $v_{\rm min}>0$ and naturally also the time horizon to be positive. In order to avoid complications coming from streets that are too short to support communications, in the first step, let us suppose that $\rho=0$ and prove existence of percolation for large device intensities. For this, let us present a condition for the waypoint kernel that relaxes isotropy around the starting position. 
We say that a kernel $\k$ is {\em $c$-well behaved} if for almost-all $S$ and all $x\in S$,
\begin{align*}
B_c(x) \cap S \subset \text{supp}\big(\k^S(x, \d y)\big).
\end{align*}
Further, we say that $\k$ is {\em well behaved} if $\k$ is $c$-well behaved for some $c > 0$.
Let us note that both $\k'$ and $\k''$ are $L$-well behaved.
%Let us note that both $\k'$ and $\k''$ are well behaved as long as $R>0$. 
Here is our first result for the immediate connectivity $\rho=0$. 
\begin{theorem}\label{thm_largeintense}
Let $\k$ be well behaved and $T>0$ then, $\lac(T,\muv,0,r)<\infty$. 
\end{theorem}
%Let us mention that this statement is in fact true under the assumption that $\muv\neq\delta_0$, as can be seen from the proof. 

The proof rests on a coupling of the original process to an auxiliary process of open streets such that existence of an unbounded component of open streets implies percolation of the original model. Here, we say that a street is \textit{open} if devices on that street can connect to devices on all its neighboring streets and also inter-connect on that street. In particular, any infinite connected component of open streets contains an infinite connected component of devices. 
However, as $\la$ increases, the probability for a street to be open converges to one, resulting in a percolation regime for open streets. 
We present the complete proof in Section~\ref{Sec_Proof_Sup}. 

\medskip
We have now exhibited the existence of a percolation regime for $\rho = 0$, if $\lambda$ is sufficiently large. In order to move beyond this setting, note that for $\rho>0$, very short streets $s$ with $|s|<\rho v_{\rm min}/2$ can never be used to establish a connection, however, they can be used for traveling. 
As a consequence, for the case where $\rho > 0$ we need stronger assumptions on $\muv$ and $r$ in order to guarantee that $\lac(T,\muv,\rho,r)<\infty$.

Again, we want to compare our original model to the simpler model on street systems $S^a$, where all streets in $S$ that have length less than $a$ are eliminated. Recall that we say that $S^a$ percolates if it contains an infinite connected component of streets and 
$$
a_{\rm c}(0)=\sup\{a>0\colon \P(S^{a}\text{ percolates})>0\}. 
$$
In principle, $a_{\rm c}(0)$ could very well be zero, however, in Theorem~\ref{thm_twicethinstreet} below, we show that indeed, $a_{\rm c}(0)>0$ for all asymptotically essentially connected street systems. Let us note that $a_{\rm c}(0)\le a_{\rm c}$, where $a_{\rm c}$ is defined below Definition~\eqref{def_abc}. In general it is not clear whether we also have that $a_{\rm c}(0)\ge a_{\rm c}$. This is a question of continuity of critical values, which is not always guaranteed. 

Unfortunately, asymptotically essentially connectedness can be easily violated by thinned street systems, since, with positive probability, there are disconnected components.
Let us give a different version of percolation that is beneficial for thinned street systems.
We denote by $R^a_x$ the distance of the furthest point from $x$ that is reachable without crossing $S^a$.
Now, if $\lim_{n\uparrow \infty}\P\big(\sup_{x \in Q_n\cap \Q^2}R^a_x<n\big) =1$ we say that the thinned graph is $R^a$-connected and define
$$a_{\rm c}^+:=\sup\{a>0\colon S^a \text{ is } R^a \text{-connected}\}.$$

With these definition, we are in the position to state the following result.

\begin{theorem}\label{thm-percolation-rho-positive}
Let $\k$ be $c$-well behaved for some $c>0$, $r,\rho, a_{\rm c}(0)>0$ and $T>2\rho$. Then, for all sufficiently small $v_{\rm min}$ we have $\lac(T,\muv,\rho,r)< \infty$. If furthermore $a_{\rm c}^+>0$, then the first statement holds for all $v_{\rm min} < \min(r\rho^{-1}, a_{\rm c}^+ \rho^{-1},c(2 \rho)^{-1})$.
\end{theorem}

Similar to the proof of Theorem~\ref{thm_largeintense}, the main argument is a comparison of our model to a street system in which we eliminate short streets, large streets and additionally perform an independent random thinning of street. The removal of short streets becomes necessary, since for positive $\rho$, short streets can never be used to establish connections. The independent thinning comes from the fact that, even for large device intensities, a street has good intrinsic connectivity properties only with high probability. In order to get uniform control on this probability, we have to disregard long streets. We present the complete proof in Section~\ref{Sec_Proof_Sup_2}. 

\medskip
The previous theorem heavily rests on the assumption that $S^{a}$ percolates. In the following theorem we will show that, indeed, there are regimes where this is the case. 
% \begin{theorem}\label{thm_twicethinstreet}
% Let $\mathcal{S}$ be a street system, then $a_{\rm c}(0)>0$. For street systems of Poisson--Voronoi or Poisson--Delaunay type we even have that $a_{\rm c}^+>0$.
% \end{theorem}
\begin{theorem}\label{thm_twicethinstreet}
Let $\mathcal{S}$ be a street system, then $a_{\rm c}(0)>0$ and $a_{\rm c}^+>0$.
\end{theorem}
Note that this can be seen as a complementary result to Theorem~\ref{Thm_Rho_2}, in which $a_{\rm c}<\infty$ is verified under general assumptions.

The proof rests on a comparison with a finitely-dependent site-percolation process on $\Z^2$, leveraging stabilization and asymptotic-essential-connectedness arguments. We present the complete proof in Section~\ref{Sec_Proof_Sup_3}.

\subsection{Simulations}\label{sec_sim}
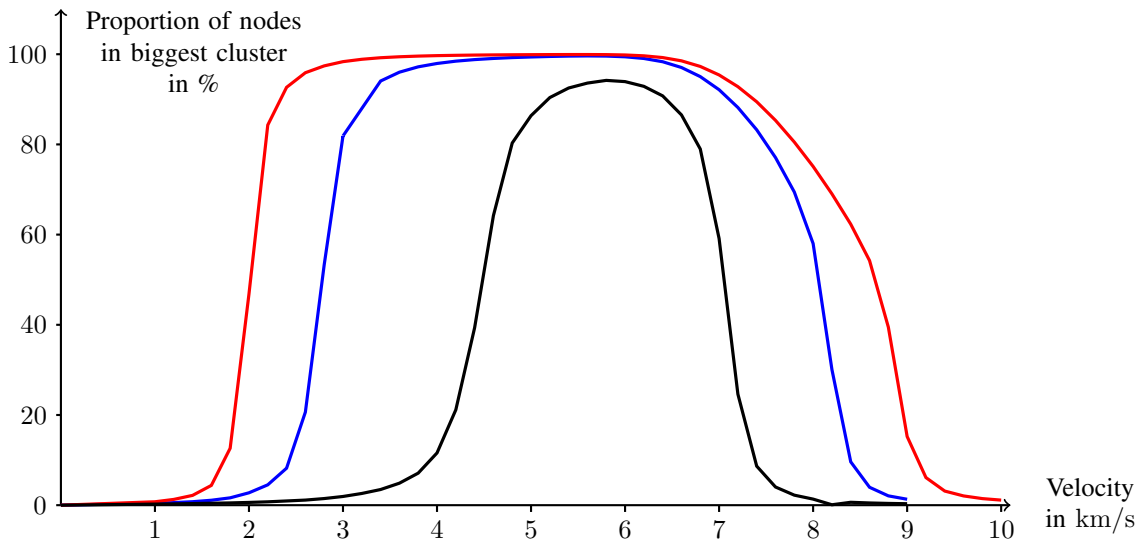
\begin{figure}[!htpb]
\centering
\begin{tikzpicture}[xscale=1.25, yscale =0.06] 
\draw[->, black, thick](0,0)--(0,110);
\draw[->, black, thick](0,0)--(10.1,0);

\coordinate[label = 0 : \begin{tabular}{c} Proportion of nodes\\  in biggest cluster \\   in \% \end{tabular}](a) at  (0.0,100);
\coordinate[label = 0 : \begin{tabular}{c} Velocity\\ in $\mathrm{km/s}$ \end{tabular}](a) at  (10.2,0);

\foreach \x in {1,2,3,4,5,6,7,8,9,10}
	\coordinate[label = -90 : $\x$](a) at  (\x,-1.5); 
\foreach \x in  {1,2,3,4,5,6,7,8,9,10}
	\draw[thick] (\x,0) -- (\x,-1.5);

\foreach \y in {0, 20,40,60,80,100}
	\draw[thick] (-0.05,\y) -- (0,\y);
\foreach \y in {0, 20,40,60,80,100}
	\coordinate[label = 180 : $\y$](a) at  (-0.05,\y);

%T=4.5
\draw[blue, very thick] plot coordinates{(0,0)(1,0.4276418921636459)(1.2,0.5711092665918871)(1.4,0.7684920031345451)(1.5999999999999999,1.0946063968287116)(1.7999999999999998,1.627970941767937)(1.9999999999999998,2.7462180815036548)(2.1999999999999997,4.513979076168805)(2.4,8.205248481173786)(2.6,20.643176904458592)(2.8000000000000003,53.281431394489466)(3.0000000000000004,81.8657208702118)(3,81.86485258175846)(3.4000000000000004,94.05927450134216)(3.6000000000000005,96.01770503156058)(3.8000000000000007,97.204584438791)(4.000000000000001,97.94775206459506)(4.200000000000001,98.46209583242849)(4.400000000000001,98.8095165437525)(4.600000000000001,99.06315741171696)(4.800000000000002,99.2524228994171)(5.000000000000002,99.38799081396265)(5.200000000000002,99.50330931851492)(5.400000000000002,99.58470963967018)(5.600000000000002,99.62250131425371)(5.8000000000000025,99.58793657143677)(6.000000000000003,99.43122938483175)(6.200000000000003,99.05640049939431)(6.400000000000003,98.3350010486429)(6.600000000000003,97.05627411836127)(6.800000000000003,95.03279577374988)(7.0000000000000036,92.1293590355353)(7.200000000000004,88.18373739186214)(7.400000000000004,83.23974464736963)(7.600000000000004,77.08008845715322)(7.800000000000004,69.44520456142371)(8.000000000000004,58.01656542635034)(8.200000000000003,30.08083548150722)(8.400000000000002,9.613587842605998)(8.600000000000001,3.970659644114552)(8.8,2.0650738945368703)(9.0,1.3136483380140012)};

%T=6
\draw [red, very thick] plot coordinates{(0,0)(1,0.7564650268799976)(1.2,1.2779744580592334)(1.4,2.1713790429056328)(1.5999999999999999,4.38540895652501)(1.7999999999999998,12.620564109570974)(1.9999999999999998,46.818073483208615)(2.1999999999999997,84.29453202377448)(2.4,92.6606516681212)(2.6,95.9006471335929)(2.8000000000000003,97.41448679796676)(3,98.35341899426601)(3.2,98.87779727947853)(3.4000000000000004,99.22334710419645)(3.6000000000000005,99.43591722477422)(3.8000000000000007,99.58436226564966)(4.000000000000001,99.68913365186794)(4.200000000000001,99.75842164936758)(4.400000000000001,99.81325238499181)(4.600000000000001,99.8495969942187)(4.800000000000002,99.88061641552293)(5.000000000000002,99.902206105413)(5.200000000000002,99.92063012845463)(5.400000000000002,99.93041659585997)(5.600000000000002,99.93176961829602)(5.8000000000000025,99.91262137952223)(6.000000000000003,99.83510906583427)(6.200000000000003,99.64670195432711)(6.400000000000003,99.26117439071664)(6.600000000000003,98.53164283631406)(6.800000000000003,97.29098066858257)(7.0000000000000036,95.40228806787293)(7.200000000000004,92.80495676458245)(7.400000000000004,89.45616652099189)(7.600000000000004,85.3298057414369)(7.800000000000004,80.52271361575305)(8.000000000000004,75.11599970311633)(8.200000000000003,69.0319490566001)(8.400000000000002,62.31539388094444)(8.600000000000001,54.24529185985242)(8.8,39.48510312098676)(9.0,15.22405800177105)(9.2,6.117704164409235)(9.399999999999999,3.129625918691648)(9.599999999999998,2.0527292606505033)(9.799999999999997,1.4584474811764477)(9.999999999999996,1.1229215785725477)};

%T=3
\draw [black, very thick] plot coordinates{(0,0)(1,0.23927743288006706)(1.2,0.2939594355176267)(1.4,0.3571183215590821)(1.5999999999999999,0.4285325643965474)(1.7999999999999998,0.5022543369636251)(1.9999999999999998,0.6128767158361121)(2.1999999999999997,0.7454406387717187)(2.4,0.9357631251031417)(2.6,1.130153111483805)(2.8000000000000003,1.4697604558399748)(3,1.930753007404037)(3.2,2.5976227331530603)(3.4000000000000004,3.4719695816811044)(3.6000000000000005,4.877020923733652)
    (3.8000000000000007,7.109181202858432)(4.000000000000001,11.584410981868079)(4.200000000000001,21.147681112466515)
    (4.400000000000001,39.35139324528766)(4.600000000000001,64.22329960862893)(4.800000000000002,80.3140540670274)
    (5.000000000000002,86.37287946674868)(5.200000000000002,90.39455201830846)(5.400000000000002,92.50979115660245)
    (5.600000000000002,93.60352791420575)(5.8000000000000025,94.19845574302252)(6.000000000000003,93.92746550842632)
    (6.200000000000003,92.87777443153484)(6.400000000000003,90.76188134677139)(6.600000000000003,86.53891372064134)
    (6.800000000000003,78.97826996558155)(7.0000000000000036,59.12755829742499)(7.200000000000004,24.615980435350654)
    (7.400000000000004,8.653208499372632)(7.600000000000004,4.029889150743293)(7.800000000000004,2.222963347764122)
    (8.000000000000004,1.345483194923647)(8.200000000000003,0.09342680221910406)(8.400000000000002,0.6422666239839645)
    (8.600000000000001,0.4930232216990354)(8.8,0.3908050746249553)(9.0,0.32287216771493894)};

\end{tikzpicture}
\caption{In-and-out of percolation: Simulation of the largest connected cluster at different times $T=3 \mathrm{min}$ (black), $T=4.5 \mathrm{min}$ (blue), and $T=6\mathrm{min}$ (red). Furthermore, we used $\rho = 10 \mathrm{sec}, r= 20\mathrm{m}$, $\lambda = 20\mathrm{devices/km}$ and a street intensity of $20\mathrm{km/km^2}$. For the velocities we chose a normal distribution, conditioned to be positive $\muv= \mathcal{N}^+(v,v/5)$. }
\label{Fig_in_out}
\end{figure}

Finally, let us present some simulations for two of the key results presented above. More precisely, we focus on the interplay between the Theorems~\ref{thm_largespeeds}, \ref{thm-percolation-rho-positive} and \ref{thm_lowintense}.
Recall that Theorem~\ref{thm_largespeeds} features absence of percolation for large velocities, which is roughly due to the fact that in this regime only very long streets can be used to form connections. On the other hand, Theorem~\ref{thm-percolation-rho-positive} establishes percolation for intermediate velocities for sufficiently large intensities of devices.
But then again, Theorem~\ref{thm_lowintense} features absence of percolation if the velocities are too small, which is due to the fact particles have insufficient time to move far enough to establish connections. 
Combining those results, we observe a double phase transition in the sense of an in-and-out of percolation in certain configurations of the parameters.
In Figure~\ref{Fig_in_out} we present a simulated graph that illustrates this behavior in a situation in which all parameters are fixed except the velocity.

Let us give some details about the simulation setting. We use PVT with an edge intensity of $20\mathrm{km/km^2}$ in order to simulate the street system of a city center, see also~\cite{courtat2012promenade} for a statistical justification of this parameter choice. We mention that, in order to reduce boundary effects, we simulate on a torus. We also note that we could effectively use Item (2) of Proposition~\ref{prop_scale} in order to significantly reduce the simulation time: As a reminder, $g_{\k,T/a,\mathcal{S},\muv^a,\rho/a,r}(X^{\la})$ percolates if and only if $g_{\k,T,\mathcal{S},\muv,\rho,r}(X^\la)$ percolates as the connections in both graphs are the same. 
Now, as
$$g_{\k,T,\mathcal{S},\muv^a,\rho,r}(X^{\la})= g_{\k,aT,\mathcal{S},\muv,a\rho,r}(X^\la),$$
 by adjusting $T$ and $\rho$, we can use the same movement from the simulation with $\muv$ if we save for each $a$ the connection graph individually. This significantly reduces the simulation load.

%
%
%
%%====================================================
%%====================================================
%
\section{Proofs}\label{sec_proofs}
%In this section we present detailed proofs of the theorems that have not been proved already. 
\subsection{Proof of Proposition~\ref{prop_scale}}\label{proof_prop_scale}

\begin{proof}[Proof of Proposition~\ref{prop_scale}]
For the spatial rescaling, note that, if we scale the initial device positions (and the underlying street system) by a factor of $a$, then we see the same waypoint selection, if we also scale the waypoint kernels such that device $aX_i$ chooses position $aY_i$ on the street system $aS$ with the same probability as $X_i$ chooses $Y_i$ on $S$, leading to the kernel $\k[a]$ and $\mathcal S^{1/a}$.
%In the street system $aS$ with $S \sim \mathcal{S}^1$ you see a factor of $a^{-2}$ less streets compared to $S$, which, however are longer by a factor of $a$. Therefore, the intensity of streets is given by $a^{-1}$.
The associated device densities per street is given by $\la/a$. Now, in order to see the same contact times, we also have to rescale the vicinity parameter $r$ by a factor $a$, and also rescale the velocity by $a$. In particular, the time horizon as well as the connection times are not rescaled, since $a$ times the length traveled with $a$ times the speed leads to the same travel times. 

On the other hand, for the time rescaling, note that, keeping the device positions, streets and waypoint kernels fixed, while rescaling the velocity by $a$, reproduces the percolation picture if we adjust the time horizon by a factor of $1/a$ and also the connection-time parameter accordingly by $1/a$. 
\end{proof}

\subsection{Proof of Theorem~\ref{thm_lowintense}}\label{Sec_Proof_Sub_Rho}
\begin{proof}[Proof of Theorem~\ref{thm_lowintense}]
Let us first prove $\lac(T,\muv,\rho,r)>0$.
Each device can travel at most to a distance $Tv_{\rm max}$ from its starting position. Any connected component of $g_{T,\muv,\rho,r}(X^\la)$ is contained in a connected component of the Cox--Gilbert graph $g_{2(T-\rho)v_{\rm max}+r}(X^\la)$, where any two devices in $X^\la$ are connected if their mutual distance is less than $2(T-\rho)v_{\rm max}+r$. 
%In the Gilbert graph, any pair of points is connected if and only if their mutual distance is smaller than ${2(T-\rho)v_{\rm max}+r}$. 
Indeed, as the devices need time $\rho$ to form a connection, the connection attempt must start before time $T-\rho$ with distance at most $r$ from each other.
 But $g_{2(T-\rho)v_{\rm max}+r}(X^\la)$ percolates with probability zero for all sufficiently small  $\la$, see, for example~\cite{wias2}. %The same argument only applies for $T$ and $v_{\rm max}$ for $r = 0$ (and $\rho = 0$ as a consequence, \red{see Theorem~\ref{Thm_Rho_0}}). For $r > 0$, it is possible that the Gilbert graph already $g_{r}(X^\la)$ percolates. Therefore, other arguments that rely on the street system are required.
 
\medskip
In order to prove $\vc(T,\la,\rho,r)>0$, note that the previous comparison argument can not be used if the Cox--Gilbert graph $g_r(X^\lambda)$ already percolates. However, using that communication on different streets is not possible, in order to see a connection between devices that have started on different streets, at least one device has to move across a crossing. Using this, we construct a coupling to a site-percolation model on $\Z^2$. 
For this, we call a crossing $C_i$ {\em speed-open} if there is at least one device within distance $\vmax T$ along the adjacent streets, and {\em speed-closed} otherwise. Note that the probability for a crossing to be speed-open is given by $1-\exp(-\lambda| B^S_{\vmax T}(C_i)\cap S|)$, where, for $x \in S$, we denote by $B^S_l(x)\subset S$ all the points within distance $l$ measured along the street system. %Although this event is not independent, it goes to zero as $\vmax \downarrow 0$. There is a $\vmax > 0$ such that the site-percolation model does not percolate. This implies non-percolation of the original model and $\vc > 0$.
Unfortunately, $|B^S_{\vmax T}(C_i)\cap S|$ is not bounded. Therefore, more refined arguments using the stationarity and regularity of the street system are required.
We are going to construct boxes where we both control the number of crossings and the number of streets per area. 
Recall the random field of stabilization radii $\{R_x\}_{x\in\R^2}$ for the segment process $S$ and the notation $R(V)=\sup_{y\in V\cap \Q^2}R_y$. 
For fixed parameters $s,m$, and $n$, we say that $z \in \Z^2$ is \textit{closed}, if we have that 
\begin{itemize}
\item[(a)] $R(Q_{6n}(nz))< n$,
\item[(b)] $|\{i \colon C_i \in Q_{6n}(nz)\}|<m$
%item[(c)] $\max_{\{C_i\}_{i \in I}\cap Q_{6n}}|B_{1}(C_i)\cap S| < s$ and
\item[(c)] $\min_{C_i\sim C_j \in Q_{6n}(nz)} |C_i - C_j| > s$
\item[(d)] every crossing in $Q_{3n}(nz)$ is speed-closed,
\end{itemize}
and {\em open} otherwise. Here we write $C_i\sim C_j$ to indicate that the crossings are associated to the same street. We denote by $A,B,C,$ and $D$ respectively the events that the conditions (a), $\dots$, (d) are not fulfilled for the origin. Then
\begin{align*}
    \P(o \mathrm{\ is \ open})=\P(A \cup B \cup C\cup D ) \le \P(A)+\P(B \cup C) + \P(D \cap A^{\mathrm{c}} \cap B^{\mathrm{c}} \cap C^{\mathrm{c}} ).
\end{align*}
As the model is stabilizing, we are able to choose $n=n(\e)$ big enough such that $P(A)<\e/3$.
Now, since $S$ is locally finite, we can choose $m=m(n,\e)$ and $s=s(n,\e)$ such that $\P(B \cup C)<\e/3$.
Finally, under the event $A^{\mathrm{c}} \cap B^{\mathrm{c}} \cap C^{\mathrm{c}}$, for $\vmax$ sufficiently small such that $\vmax T <s$ and $C_i \in Q_{3n}$ we can bound $$|B^S_{\vmax T}(C_i)| = \mathrm{deg}(C_i) \vmax T \le m \vmax T.$$ For the equality we used that, due to (c), no other crossing is within distance $\vmax T$. Due to (a), $C_i\in Q_{3n}(nz)$ is not connected to a crossing outside $Q_{6n}(nz)$. Therefore $\mathrm{deg}(C_i)$ can be bounded by $m$ due to (b). For the special case where $S$ is a PVT this much care is not required as $\mathrm{deg}(C_i) = 3$. Now,
\begin{align*}
    \P(D \cap A^{\mathrm{c}} \cap B^{\mathrm{c}} \cap C^{\mathrm{c}}) 
    < m (1-\exp(-\lambda m \vmax T)).
\end{align*}
Therefore, we can choose $\vmax=\vmax(n,m,s, \epsilon)$ small enough such that 
$\P(D \cap A^{\mathrm{c}} \cap B^{\mathrm{c}} \cap C^{\mathrm{c}})<\e/3$.
Therefore, 
\begin{align} \label{eq_proof_largeintense}
\limsup_{n\uparrow\infty}\limsup_{m\uparrow\infty}\limsup_{s\downarrow0}\limsup_{\vmax \downarrow0}\P( o \text{ is open} )=0.
\end{align}
Now, in order to see absence of percolation of open sites in $\Z^2$, note that the process of open sites, due to condition (a) (which implies stabilization), is $6n$-dependent and stationary due to the translation covariance of the waypoint kernel and the street system. Thus, we can invoke the domination-by-product-measures theorem \cite[Theorem 0.0]{domProd} to dominate the process of open sites, for sufficiently small $\vmax$, by a subcritical Bernoulli percolation process on $\Z^2$. 
Finally, note that absence of percolation of open sites implies absence of percolation in our original model.

\medskip
In order to prove $\Tc(\la,\muv,\rho,r)>2\rho$, we compare to a similar site-percolation model.
As the original model is monotone in $T$, we show existence of $\e>0$ such that for $T=2\rho +\e$ we see no percolation.
In order for a device to transmit on two different streets, it has to spend a time of $\rho$ on each street it transmits on. Therefore, any street change of the device has to occur in the time interval $[\rho, \rho+\e]$.

We say that a crossing $C_i$ is {\em time-open} if there is a device  $X_j$ such that both $T_{j,t} = C_i$ for a $t \in [\rho, \rho+\e]$  and there is an time interval of length at least $\rho$ that the device spends on a street. We call the crossing {\em time-closed} otherwise.
Similar to the arguments for $\vmax \downarrow 0$, a control on the regularity of the street system yields the desired result that site percolation is not possible if $\e$ is chosen small enough.
For this, let us replace Condition (d) above by the condition
\begin{itemize}
    \item[(d')] every crossing in $Q_{3n}(nz)$ is time-closed,
\end{itemize}
and say that $z \in \Z^2$ is \textit{closed} if (a), (b), (c) and (d') are fulfilled, and {\em open} otherwise.
It remains to prove a uniform bound to see a time-open crossing in order to prove that the open vertices in $\Z^2$ do not percolate.
In order to ensure that no device moves across the crossing in a fixed time-interval, we have to ensure that no device is positioned on certain intervals on the streets. Unfortunately, the position of the interval depends on the velocities of the devices. For devices of fixed speed $\vmin \le v\le \vmax$, we can ensure that no device moves across the crossing $C_i$ if there are no devices positioned in $B^S_{v(\rho+\e)}(C_i)\setminus B^S_{v\rho}(C_i)$. Here, we can ignore devices $X_j \in B^S_{v\rho}(C_i)$, as any sub-path of a shortest path is a shortest path itself and therefore can not spend a time of $\rho$ on any street before crossing $C_i$, almost surely.

Let us assume that $\e$  is small enough such that $\e \vmax < s$. As $\e\vmax$ is lower than the minimum street length, it is not possible for a device to fully traverse one street of length $v\rho$ and (before or afterwards) fully traverse an additional street.
As there is at most one crossing, using $\mathrm{deg}(C_i) < m$ for all $C_i \in Q_{3n}(nz)$ we obtain for all speeds $\vmin \le v \le \vmax$ the uniform bound 
$$|B^S_{v(\rho+\e)}(C_i)\setminus B^S_{v\rho}(C_i)|<m^2\vmax\e,$$
where we used $\deg(C_i)<m$ to bound the number of streets.
We then sum over all the crossings in $Q_{3n}\le m$. Interchanging the order of integration of the Cox point process $X^\la$ and its attached paths, Fubini's Theorem yields that
\begin{align*}
    \P(D' \cap A^{\mathrm{c}} \cap B^{\mathrm{c}} \cap C^{\mathrm{c}})  & \le m \E\Big [  1-\exp\big(-\int\muv(\d v)\lambda|B^S_{v(\rho+\e)}(C_i)\setminus B^S_{v\rho}(C_i)|\big) \one\{ A^{\mathrm{c}} \cap B^{\mathrm{c}} \cap C^{\mathrm{c}}\}\Big] \\
    &\le m (1-\exp(-\lambda m^2 \vmax \e)) ,
\end{align*}
which converges to $0$ as $\e \downarrow 0$.

Using the same domination argument as above, sufficiently small $\e$ leads to absence of percolation of good sites. Further, if the site-percolation process of open sites does not percolate, there is also no percolation in the original model.
Choosing $\e>0$ small enough that the site-percolation process is subcritical, implies $\Tc(\la,\muv,\rho,r)\ge 2\rho+\epsilon>2\rho$.
% As the connection model is monotone in $\rho$, e.g. $g_{\k,T,\g,\muv,\rho,r}(X^\la) \subset g_{\k,T,\g,\muv,0,r}(X^\la)$
% Due to the temporal scaling relation (item (2) of Proposition~\ref{prop_scale}, we obtain 
% $g_{\k,T/a,\g,\muv^a,\rho/a,r}(X^{\la})$ percolates if and only if $g_{\k,T,\g,v,\rho,r}(X^\la)$ percolates.
% \end{enumerate}
\end{proof}

\subsection{Proof of Theorem~\ref{thm_largespeeds}}\label{Sec_Proof_Sub}
For $x,y\in \R^2$ denote by $\ell_S(x,y)\subset \R^2$ the shortest path between $x$ and $y$ on $S$. Further denote by $\ell'_S(x)=\sup\{|\ell_S(x,y)|\colon y\in\supp(\k^S(x,\d y))\}$ the length of the shortest path starting in $x$ towards any reachable target on $S$. It will be convenient to write $\ell_S(x)=\ell'_S(x)+r/2$ for the extended shortest-path lengths. Further, we consider the following set
\begin{align*}
A'_S(v)=\{x\in S\colon &\text{ there exists }y\in\supp(\k^S(x,\d y))\text{ such that }\\ 
&\ell_S(x,y)\text{ contains a street }s\text{ of length }|s|>v\rho/2\},
\end{align*}
and note that it only depends on the street system $S$, but not on the point process of devices. 
Denote by 
$$X^{v,S,\la}=\{X_i\in X^\la\colon X_i\in A'_S(v)\}$$ 
the thinned version of $X^\la$, where the thinning rule depends only on $S$ and $v$.
Then, we consider the following {\em geostatistical Cox--Gilbert graph} $h_{v}(X^\la)$ with vertex set $X^{v,S,\la}$ and edges between any pair of vertices $X_i,X_j\in X^{v,S,\la}$ if and only if 
$$|X_i-X_j|\le\ell_S(X_i)+\ell_S(X_j).$$
Note that indeed, the interaction range $\ell_S(X_i)$ associated to $X_i$ is determined by the underlying street system. Additionally, the thinning of the Cox point process is also determined by the street system. 

Any edge in $g_{v,\rho,r}(X^\la)$ is also contained in $h_{v}(X^\la)$ and thus it suffices to prove absence of percolation in $h_{v}(X^\la)$. 

Our approach rests on a generalization of arguments first presented in~\cite{G08}, for Poisson--Boolean models and their extensions towards Cox--Boolean models in~\cite{JTC20}, which leverage scaling properties of the process and adjustments towards our geostatistical setting. It will be thus convenient to reformulate our results in terms of the associated {\em geostatistical Cox--Boolean model}
$$\C=\bigcup_{X_i\in X^{v,S,\la}}B_{\ell_S(X_i)}(X_i),$$
where we note that $\C$ contains an infinite connected component, if and only if $h_v(X^\la)$ contains an infinite connected component. 

Let $\C_x(V)$ denote the connected component containing $x$ of the geostatistical Cox--Boolean model, based on points in $V\subset\R^2$. Then, we define for any $x\in \R^2$ and $\alpha>0$ the event
\begin{align*}
G(x,\a)&=\big\{\C_x(B_{10\a}(x))\not\subset B_{8\a}(x)\big\}
\end{align*}
that the cluster of $x$, only using points in $B_{10\a}(x)$, reaches beyond $B_{8\a}(x)$. Consider the event
\begin{align*}
A_S(\a)=\{x\in S\colon \ell_S(x)<\a \},
\end{align*}
and note that $A_S$ only depends on the street system $S$ and $\k$. 
Let $S^*$ denote the Palm version of $S$, and recall the stabilization radii $R_y$ of the street system. We write $R(V)=\sup_{y\in V\cap\, \Q^2}R_y$ for any $V\subset\R^2$. Then, we have the following key lemma that establishes a scaling relation in the model.
\begin{lemma}\label{Lem_A}
Consider the geostatistical Cox--Boolean model. Then, there exists a constant $c>0$ such that for all $\a>0$, we have 
\begin{align*}
\P\big(G(o,10\a)\big)&\le c\P \big(G(o,\a)\big)^2+c\a^2\P\big(o\in A^{\rm c}_{S^*}(\a)\big)+c\P\big(R(Q_{10\a})\ge \a\big)
\end{align*}
\end{lemma}
\begin{proof}[Proof of Lemma~\ref{Lem_A}]
For a measurable set $B\subset \R^2$, consider the event
\begin{align*}
F_B(\a)=\{&\text{for all }X_i\in X^\la\cap B\text{ we have that }X_i\in A_S(\a) \},
\end{align*}
and let  $K_\a$ denote a finite subset of the sphere $\S_\a=\{x\in\R^d\colon |x|=\a\}$ such that 
$$\S_\a\subset \bigcup_{x\in K_\a}B_1(x).$$
Then, using arguments first presented in~\cite{G08,JTC20} gives 
\begin{align*}
\P\big(G(o,10\a)\big)&\le\P\Big\{\Big( \bigcup_{k\in K_{10}} \big(G(\a k,\a)\cap F_{B_{10\a}(\a k)}(\a)\big)\Big) \cap \Big(\bigcup_{l\in K_{50}}\big(G(\a l,\a)\cap F_{B_{10\a}(\a l)}(\a)\big)\Big)\Big\}\\
&\qquad+\P(F_{B_{100\a}}(\a)^{\rm c})\\
&\le\sum_{k\in K_{10},\, l\in K_{50}}\E\big[\P\big(G(\a k,\a)\cap F_{B_{10\a}(\a k)}(\a)|S\big)\P\big(G(\a l,\a)\cap F_{B_{10\a}(\a l)}(\a)|S\big)\big]\\
&\qquad+\P(F_{B_{100\a}}(\a)^{\rm c}),
\end{align*}
since $G(o,10\a)\subset \Big(\Big(\bigcup_{k\in K_{10}} G(\a k,\a)\Big)\cap \Big(\bigcup_{l\in K_{50}}G(\a l,\a)\Big)\Big)$ and under the event $F_{B_{100\a}}(\a)$, the events are independent, conditioned on $S$. Indeed, if a path exists with length smaller than $\a$, then also the shortest path must be contained in a disc of radius $\a$ around the starting position. The idea here is that if the process can percolate beyond $B_{80\a}(o)$ and the individual traveling ranges are all small, then the process must pass through two circles of discs and percolate in at least one smaller disc in the inner circle and in the outer circle. Introducing the stabilization events, we can further bound this and arrive at the desired expression using also stationarity, once we realized that, by Campbell's formula,
\begin{align*}
\P(F_{B_{100\a}}(\a)^{\rm c})&\le \E\Big[\sum_{X_i\in X^\la\cap B_{100\a}}\one\{X_i\in A^{\rm c}_S(\a)\}\Big]\\
&\le \la\E\Big[\int_{S\cap B_{100\a}}\d x\one\{x\in A^{\rm c}_S(\a)\}\Big]=\la\pi 100^2\a ^2\P(o\in A^{\rm c}_{S^*}(\a)),
\end{align*}
where we used that the street system is normalized.
\end{proof}
Note that, since we assume stabilization, we have that $\P(R(Q_{10\a})\ge \a)$ tends to zero as $\a$ tends to infinity. The same is true for the other error term, as we show now. 
\begin{lemma}\label{Lem_B}
We have that $ \a ^2\P\big(o\in A^{\rm c}_{S^*}(\a)\big)$ tends to zero as $\a$ tends to infinity. 
\end{lemma}
\begin{proof}[Proof of Lemma~\ref{Lem_B}]
The idea is to first ensure existence of paths via asymptotic essential connectedness. Then, we bound the length of that path via the total street length. This is achieved by coupling the number of streets to the number of faces with the Euler formula. Using different scales for $\a$, we arrive at the desired result. 

First, we have that $\P\big(o\in A^{\rm c}_{S^*}(\a)\big)=\a^{-1/2}\E\Big[\int_{S\cap Q_{\a^{1/4}}}\d x\one\{x\in A^{\rm c}_S(\a)\}\Big]$. Next, 
\begin{align*}
\E\Big[\int_{S\cap Q_{\a^{1/4}}}\d x\one\{x\in A^{\rm c}_S(\a)\}\Big]&=\E\Big[\int_{S\cap Q_{\a^{1/4}}}\d x\one\{x\in A^{\rm c}_S(\a)\}\one\{R(Q_{4\a^{1/4}})<\a^{1/4}\}\Big]+ O\big(\exp(-\alpha^{1/4})\big),
\end{align*}
% \red{$v-e+f = 2$
% $3v=2e$}
where the error term vanishes exponentially fast in $\a^{1/4}$ since we work with PVT or PDT, see~\cite{wias2} for details.  For the main term, under the event for the stabilization radii and sufficiently large $\a$ (larger than the support of $\k$), for every pair of points in $x,y\in S\cap Q_{2\a^{1/4}}$ there exists a path within $S\cap Q_{4\a^{1/4}}$ and this path has maximal length given by $|S\cap B_{4\a^{1/4}}|$. This is due to the asymptotic-essential-connectedness property of the PVT and PDT. Hence we can bound, 
\begin{align*}
\E\Big[\int_{S\cap Q_{\a^{1/4}}}&\d x\one\{x\in A^{\rm c}_S(\a)\}\one\{R(Q_{4\a^{1/4}})<\a^{1/4}\}\Big]\\
&\le\E\Big[|S\cap Q_{\a^{1/4}}|\one\{|S\cap B_{4\a^{1/4}}|>\a\}\one\{R(Q_{2\a^{1/4}})<\a^{1/4}\}\Big]. 
\end{align*}
Since PVT is a simple planar graph where every vertex has degree 3, we can employ Euler's formula and bound the number of streets in $S\cap Q_{\a^{1/4}}$ by $3Y(Q_{2\a^{1/4}})$ and obtain
\begin{align*}
|S\cap Q_{\a^{1/4}}|\le 6\a^{1/4}Y(Q_{2\a^{1/4}}),
\end{align*}
where we used that every street in $Q_{\a^{1/4}}$ has maximal length given by $\sqrt{2}\a^{1/4}$ and $Y$ denotes the underlying Poisson point processes for the street system, in particular $Y(A)$ denotes the number of points of $Y$ in $A\subset \R^2$, i.e., is a Poisson random variable. Note that we also used the stabilization event here.  Hence, introducing another stabilization event, we can further bound, 
\begin{align*}
\E\Big[\int_{S\cap Q_{\a^{1/4}}}&\d x\one\{x\in A^{\rm c}_S(\a)\}\one\{R(Q_{4\a^{1/4}})<\a^{1/4}\}\Big]\\
&\le6\a^{1/4}\E\Big[Y(Q_{2\a^{1/4}})\one\{Y(Q_{8\a^{1/4}})>\frac{1}{24}\a^{3/4}\}\Big]+O(\exp(-\a^{1/4}))\\
&\le 48\a^{3/4}\P(Y(Q_{8\a^{1/4}})>\frac{1}{24}\a^{3/4}-1)+O(\exp(-\a^{1/4})). 
\end{align*}
Using Poisson concentration inequalities, we see that the last probability tends to zero exponentially fast in $\a^{1/4}$, which gives the desired convergence. In the case where PVT is replaced by PDT, we note that every face in the PDT is a triangle and hence another application of Euler's formula implies the result. 
\end{proof}
Next we show that the local percolation probability becomes zero for large velocities. 
\begin{lemma}\label{Lem_C}
There exists $c>0$ such that $\P\big(G(o,\a)\big)\le  c\a^2\P\big(o\in A'_{S^*}(v)\big)$ and $\P\big(o\in A'_{S^*}(v)\big)$ tends to zero as $v$ tends to infinity. 
\end{lemma}
\begin{proof}[Proof of Lemma~\ref{Lem_C}]
Note that 
\begin{align*}
\P\big(G(o,\a)\big)&\le \P\big(X^{v,S,\la}(B_{10\a})>0\big)\le \E\big[X^{v,S,\la}(B_{10\a})\big]\\
&\le \la\E\Big[\int_{S\cap B_{10\a}}\d x\one\{x\in A'_S(v)\}\Big]\le \la 10^2\a^2\P\big(o\in A'_{S^*}(v)\big).
\end{align*}
Using the dominated-convergence theorem, it suffices to show that $\lim_{v\uparrow\infty}\one\{o\in A'_{S^*}(v)\}=0$ for almost-all realizations of $S^*$. But this is true by the almost-sure uniqueness of shortest paths in $S^*$. 
\end{proof}
Next we show that indeed, the probability to see large connected components can be bounded by the local percolation probability and some error term. 
\begin{lemma}\label{Lem_D}
Consider the geostatistical Cox--Boolean model $\C$. Then, there  exists a constant $c>0$, such that for all $\a>0$, we have 
\begin{align*}
\P\big(X^{v,S,\la}(\C_o(\R^2))> X^{v,S,\la}(B_{8\a})\big)&\le \P\big(G(o,\a)\big)+\la\E\Big[\int_{S}\d x\one\{10\a\le |x|\le 9\a+\ell_S(x)\}\Big].
\end{align*}
\end{lemma}
\begin{proof}[Proof of Lemma~\ref{Lem_D}]
Consider the event, 
\begin{align*}
H(\a)&=\{\text{there exists }X_i\in X^{v,S,\la}\cap B^{\rm c}_{10\a}\colon |X_i|\le 9\a+\ell_S(X_i)\}.
\end{align*}
Then, we have that $\big(G(o,\a)^{\rm c}\cap H(\a)^{\rm c}\big)\subset \{\C_o(\R^2)\subset B_{8\a}\}$ since, under the event $H(\a)^{\rm c}$, points outside $B_{9\a}$ cannot help the cluster $C_o$ to reach outside of $B_{8\a}$. But since $\{\C_o(\R^2)\subset B_{8\a}\}\subset \{X^{v,S,\la}(\C_o(\R^2))\le X^{v,S,\la}(B_{8\a})\}$, we have 
\begin{align*}
\P\big(X^{v,S,\la}(\C_o(\R^2))> X^{v,S,\la}(B_{8\a})\big)&\le \P\big(G(o,\a)\big)+\P\big(H(\a)\big). 
\end{align*}
Finally, using Campbell's formula
\begin{align*}
\P\big(H(\a)\big)&\le \la\E\Big[\int_{S}\d x\one\{10\a\le |x|\le 9\a+\ell_S(x)\}\Big],
\end{align*}
which finishes the proof.
\end{proof}
The error term becomes small for large $\a$, this is the content of the next result. 
\begin{lemma}\label{Lem_E}
We have that $\E\Big[\int_{S}\d x\one\{10\a\le |x|\le 9\a+\ell_S(x)\}\Big]$ tends to zero as $\a$ tends to infinity. 
\end{lemma}
\begin{proof}[Proof of Lemma~\ref{Lem_E}]
We again want to employ asymptotic essential connectedness in order to bound $\ell_S(x)$ against the total street length in a disc. We need slightly more sophisticated arguments compared to the proof of Lemma~\ref{Lem_B}. 
We want to cover $B_{10\a}^{\rm c}$ by layers of squares with side-length of smaller order than $\a$ where each layer has a growing diameter such that also the number of squares in the layers increases polynomially.
For this, consider first a sequence of auxiliary squares $Q'_{n^4\a^{1/4}}$ and the associated partition of $\R^2$ given by the union of annuli $Q'_{(n+1)^4\a^{1/4}}\setminus Q'_{n^4\a^{1/4}}$. We call this annulus $Q'_{(n+1)^4\a^{1/4}}\setminus Q'_{n^4\a^{1/4}}$ the $n$-th layer. Then, we cover the $n$-th layer by shifted non-overlapping squares of the form $Q_{n\a^{1/4}}$ and note that this can be done with $c'n^5$ many squares for some constant $c'>0$. We write  $Q'_{(n+1)^4\a^{1/4}}\setminus Q'_{n^4\a^{1/4}}\subset \bigcup_{j\in U(n,\a)}Q^j_{n\a^{1/4}}$ for such a covering. In particular, 
\begin{align*}
\E\Big[&\int_{S}\d x\one\{10\a\le |x|\le 9\a+\ell_S(x)\}\Big]\\
&\le \sum_{n\ge 0}\sum_{j\in U(n,\a)}\E\Big[\int_{S\cap Q^j_{n\a^{1/4}}}\d x\one\{10\a\le |x|\le 9\a+\ell_S(x)\}\Big].
\end{align*}
Now we introduce the asymptotic-essential-connectedness event $\{R(Q^j_{4n\a^{1/4}})<n\a^{1/4}\}$ and its complement, leading to a sum of two term, which we estimate separately. For the term, under the asymptotic-essential-connectedness event, we can use similar estimates as in the proof of Lemma~\ref{Lem_B} via the total street length, to see that
\begin{align*}
\E\Big[&\int_{S\cap Q^j_{n\a^{1/4}}}\d x\one\{10\a\le |x|\le 9\a+\ell_S(x)\}\one\{R(Q^j_{4n\a^{1/4}})<n\a^{1/4}\}\Big]\\
&\le \E\Big[|S\cap Q_{n\a^{1/4}}|\one\{|S\cap Q_{4n\a^{1/4}}|\ge (n^4\a^{1/4}/2\vee 10\a)-9\a\}\Big],
\end{align*}
where we used that $|x|>n^4\a^{1/4}/2$ for all $x\in S\cap Q^j_{n\a^{1/4}}$. Again, using Euler's formula, we can further estimate against Poisson random variables $Y$, 
\begin{align*}
\E\Big[&|S\cap Q_{n\a^{1/4}}|\one\{|S\cap Q_{4n\a^{1/4}}|\ge (n^4\a^{1/4}/2\vee 10\a)-9\a\}\Big]\\
&\le \sqrt{2}\, 6n\a^{1/4}\E\Big[Y(Q_{8n\a^{1/4}})\one\{Y(Q_{8n\a^{1/4}})\ge \frac{(n^4\a^{1/4}\vee 10\a)-9\a}{\sqrt{2}\, 48n\a^{1/4}}\}\Big].
\end{align*}
In particular, using Poisson concentration inequalities, we see that the last expression tends to zero exponentially fast, both in $n$ and in $\a^{1/4}$. The crucial observation here is that, in the indicator on the right-hand side, the lower bound on the right increases on a larger order than the intensity of the Poisson random variable, both in $n$ as well as in $\a$. In particular, there exists constants $c_1,c_2>0$, such that,
\begin{align*}
\sum_{n\ge 0}&\sum_{j\in U(n,\a)}\E\Big[\int_{S\cap Q^j_{n\a^{1/4}}}\d x\one\{10\a\le |x|\le 9\a+\ell_S(x)\}\one\{R(Q^j_{4n\a^{1/4}})<n\a^{1/4}\}\Big]\\
&\le c_1\sum_{n\ge 0} n^7\a^{3/4}\exp(-c_2 n\a^{1/4}),
\end{align*}
which tends to zero as $\a$ tends to infinity. 

For the error term, due to absence of stabilization $\{R(Q^j_{4n\a^{1/4}})\ge n\a^{1/4}\}$, note that this term vanishes exponentially fast also both in $n$ and $\a^{1/4}$ and hence, there exists constants $c_1,c_2,c_3>0$, such that,
\begin{align*}
\sum_{n\ge 0}&\sum_{j\in U(n,\a)}\E\Big[\int_{S\cap Q^j_{n\a^{1/4}}}\d x\one\{10\a\le |x|\le 9\a+\ell_S(x)\}\one\{R(Q^j_{4n\a^{1/4}})\ge n\a^{1/4}\}\Big]\\
&\le c_1 \sum_{n\ge 0}n^6\a^{1/4}\E\big[Y(Q_{n\a^{1/4}})\one\{R(Q_{4n\a^{1/4}})\ge n\a^{1/4}\}\big]\\
&\le c_2 \sum_{n\ge 0}n^8\a^{3/4}\P^{1/2}\big(R(Q_{4n\a^{1/4}})\ge n\a^{1/4}\big)\\
&\le c_2\sum_{n\ge 0} n^8\a^{3/4}\exp(-c_3 n\a^{1/4}),
\end{align*}
where we used H\" older's inequality and the fact that the second moment of a Poisson random variable with parameter $b$ is given by $b^2+b$. 
But the last expression tends to zero as $\a$ tends to infinity, which finishes the proof.  
\end{proof}
We will also need the following essential result about convergence properties of functions satisfying some scaling inequality. 
\begin{lemma}[{\cite[Lemma 3.7]{G08}}]\label{Lem_F}
Let $f$ and $g$ be two bounded measurable functions from $[1,\infty]$ to $[0,\infty)$. Additionally, let $f$ be bounded by $1/2$ on $[1,10]$, $g$ be bounded by $1/4$ on $[1,\infty]$ and assume 
\begin{align*} 
f (\a) \le f (\a/10)^2 + g(\a),\qquad\text{for all }\a\ge 10.
\end{align*}
Then, $\lim_{\a\uparrow\infty}g(\a)=0$ implies that $\lim_{\a\uparrow\infty}f(\a)=0$. 
\end{lemma}
We have now assembled all necessary tools in order to proof existence of a subcritical percolation regime for large speeds. 
\begin{proof}[Proof of Theorem~\ref{thm_largespeeds}]
We assume that $S$ to be a PVT or PDT. In order to prove that $\vcc<\infty$, it suffices to show that $\lim_{\a\uparrow\infty}\P(X^\la(C_o)> X^\la(B_{8\a}))=0$ for all sufficiently large $v$. But this is true if $\lim_{\a\uparrow\infty}\P\big(G(0,\a)\big)=0$, by an application of Lemma~\ref{Lem_D} and Lemma~\ref{Lem_E}. In order to show $\lim_{\a\uparrow\infty}\P\big(G(0,\a)\big)=0$ we apply Lemma~\ref{Lem_A}, Lemma~\ref{Lem_B}, Lemma~\ref{Lem_C} and Lemma~\ref{Lem_F} for proper choices of $f$ and $g$. 

For this, first define 
\begin{align*}
\a_1&:=\inf\{s\ge 1\colon \phi_1(\a)=\a ^2\P\big(o\in A^{\rm c}_{S^*}\big)<(8c^2)^{-1}\text{ for all }\a\ge s\}\\
\a_2&:=\inf\{s\ge 1\colon \phi_2(\a)=\P\big(R(Q_{10\a})\ge \a\big)<(8c^2)^{-1}\text{ for all }\a\ge s\},
\end{align*}
and define $\a_{\rm c}=\a_1\vee\a_2$. Note that $\a_{\rm c}<\infty$ by the assumption of stabilization and Lemma~\ref{Lem_B}. Based on this, we make the following definitions. We define
\begin{align*}
v'^{\rm c}&=\inf\{v\ge 1\colon \P\big(o\in A'_{S^*}(v)\big)<\frac 1 2 (100c\a_{\rm c})^{-2}\text{ for all }v\ge s\},
\end{align*} 
where $\vcc<\infty$ by Lemma~\ref{Lem_C}. With this, define 
\begin{align*}
f(\a)=c\P\big(G(o,10\a_{\rm c}\a)\big)\quad&\text{ and }\quad g(\a)=c^2\big(\phi_1(\a_{\rm c}\a)+\phi_2(\a_{\rm c}\a)\big).
\end{align*}
Then, indeed, using again Lemma~\ref{Lem_C}, we have that 
\begin{align*}
f(\a)&\le \frac{1}{2},\qquad\text{ for all }1\le \a\le10\text{ and }v>v'^{\rm c}.
\end{align*}
Further, using Lemma~\ref{Lem_B}, we have that 
\begin{align*}
g(\a)&\le \frac 1 4,\qquad\text{ for all }1\le \a.
\end{align*}
Finally, using Lemma~\ref{Lem_A}, we have that 
\begin{align*}
f(\a)&\le c^2\big(\P\big(G(o,\a_{\rm c}\a)\big)^2+\phi_1(\a_{\rm c}\a)+\phi_2(\a_{\rm c}\a) \big)=f(\a/10)^2+g(\a)\quad\text{ for all }\a \ge 10.
\end{align*}
Hence, since $\lim_{\a\uparrow\infty}g(\a)=0$, an application of Lemma~\ref{Lem_F} gives the result. 
\end{proof}

\subsection{Proof of Theorem~\ref{Thm_Rho_1}}\label{Sec_Proof_Sub_7}

\begin{proof}[Proof of Theorem~\ref{Thm_Rho_1}]
The central idea is that, even in a best-case scenario, where the device intensity is extremely large and the vicinity threshold is very small, in order for devices to form connections on at least two different streets, we must have that $2\rho+b/v\le T$, where $b$ is the graph distance between two endpoints of streets that are sufficiently long, and $v$ is a speed in the support of the speed distribution $\muv$. On the other hand, if every device can form a connection only on one single street, then the connected components in the original graph can only consist of either individual devices (that are unable to connect to any other device) or precisely those (finitely-many) devices associated to the street on which they form their connection. In particular, the initial position of a device does not have to lie on the street on which it establishes its connection. 
%If there is no other long street within distance $b$ of our long street, any device that formed a connection on this long street can only be connected to other devices that formed their connections on the same street, resulting in a local cluster.
Hence, if $S^{v_{\rm min}\rho/2,v_{\rm max}(T-2\rho)}$ does not percolate, then also our system cannot percolate.
\end{proof}

\subsection{Proof of Theorem~\ref{Thm_Rho_2}}\label{Sec_Proof_Sub_2}
\begin{proof}[Proof of Theorem~\ref{Thm_Rho_2}]
First, we will prove that for any $\mathcal{S}$ we have that $a_{\rm c}<\infty$. It suffices to exhibit $0<a,b<\infty$ such that $\P(S^{a,b}\text{ percolates})=0$. For this we employ stabilization. More precisely, we construct a subcritical percolation process of finite-dependent open sites in $\Z^2$ that dominates $S^{a,b}$.
For this, recalling the stabilization radii $\{R_x\}_{x\in\R^2}$ and the notation $R(V)=\sup_{x\in V\cap \Q^2}R_x$, we say that a site $z \in \Z^2$ is \textit{$n$-open} if
\begin{itemize}
\item[(a)] $R(Q_{3n}(nz))< n$ and
\item[(b)] $Q_{3n}(nz)\cap S^{a} =\emptyset$.
\end{itemize} 
We say that $z \in \Z^2$ is \textit{$n$-closed} if it is not $n$-open. In particular, percolation of $S^{a,b}$ implies percolation of $n$-closed sites once we pick $n>b$. Then, in order to establish absence of percolation for the process of $n$-closed sites, note that the stabilization event in Condition (a) implies that the process of $n$-closed sites is $3$-dependent. Since the whole system is also translation invariant, we can use the domination-by-product-measures theorem~\cite[Theorem 0.0]{domProd} to establish a subcritical Bernoulli site percolation on $\Z^2$ once we have shown that 
\begin{align*}
\limsup_{n\uparrow\infty}\limsup_{a\uparrow0}\P( o \text{ is $n$-closed} )=0.
\end{align*}
For this, let us denote by $A(n)$ the event described in Condition (a) and by $B(n,a)$ the event described in Condition~(b), for $z=o$. Then, it suffices to show that
\begin{align*}
\limsup_{n \uparrow\infty}\P\big( A^{\rm c}(n)\big)=0,
\end{align*}
and that for all $n\in \N$,
\begin{align*}
\limsup_{a \uparrow\infty}\P\big( B^{\rm c}(n,a)\cap A(n)\big)=0. 
\end{align*}
But, the first statement is true since we assume stabilization. For the second statement, note that, under the event $A(n)$, any street $s\subset S$ such that $s\cap Q_{n}\neq\emptyset$ has maximal length given by $|s|<\sqrt{2\, }3n$ and hence for $a>\sqrt{2}\, 3n$, we have $\P\big( B^{\rm c}(n,a)\cap A(n)\big)=0$.

\medskip
It remains to prove that, under the conditions, $\rhoc'(T,\muv)< T/2$. For this, note that, under the assumption that $a_{\rm c}<\infty$, since $\vmin>4 a_{\rm c}/T$, there exist $\eps,b>0$ such that $\vmin(T/2-\eps)/2>a_{\rm c}(b)$. In particular, $\P(S^{\vmin(T/2-\eps)/2,b'}\text{ percolates})=0$ for all $b'<b$. 
The fact that $S^{a,b}$ is decreasing in $a$ and increasing in $b$ yields for $\max(T/2-\eps,T/2-b/(2\vmax))<\rho$ that $\P(S^{\vmin\rho/2,\vmax(T-2\rho)}\text{ percolates})=0$. 
Thus, since $\max(T/2-\eps,T/2-b/(2\vmax))<T/2$, we have that $\rhoc'(T,\muv)< T/2$, which finishes this part of the proof. 
\end{proof}

\subsection{Proof of Theorem~\ref{thm_largeintense}}\label{Sec_Proof_Sup}

\begin{proof}[Proof of Theorem~\ref{thm_largeintense}]
First note that it suffices to assume that $r=0$ since larger values of $r$ only increase the percolation thresholds. 
The first main step consists in constructing a coupled process of good streets that will support percolation.
For a street $s$, let us denote by $s_1,\dots,s_n$ all neighboring streets of $s$. Then, we will call $s$ \textit{open}, if we can find devices $X^s_1,\dots,X^s_n$ on $s$ at time $t = 0$ such that
\begin{itemize}
\item[(1)]$X_1^{s},\dots,X_n^s$ are connected in  $g_{T,\muv,0,0}(X^\la \cap s)$ and 
\item[(2)] for all $i\in \{1,\dots, n\}$, there are times $0<t_i< t'_i<T$ such that $\tra_{X_i^s,t} \in s$ for $t \in [0,t_i]$ and $\tra_{X_i^s,t}^s \in s_i$ for $t \in [t_i,t'_i]$.
\end{itemize}
In words, a street is open if it contains initial positions of devices that are inter-connected on the street and where each device then travels towards its associated neighboring street, see Figure~\ref{Fig_good_street} for an illustration. 
\begin{figure}[!htpb]
\centering
\begin{tikzpicture}[scale=2.0] 
\draw[red, ultra thick](-2,0)--(1,0);
\draw[red, thick](-2,0)--(-3,2);
\draw[red, thick](-2,0)--(-3,-1);
\draw[red, thick](1.5,1)--(1,0);
\draw[red, thick](2.5,-1)--(1,0);

\draw [->, bend left, thick]  (-1,0) to (-2.4,1);
\draw[fill = blue ,ultra thin] (-1,0) circle(2pt) ;

\draw [->, bend right, thick]  (-1.4,-0.0) to (-2.3,-0.5);
\draw[fill = blue ,ultra thin] (-1.4,0) circle(2pt) ;

\draw [->, bend right, thick]  (0,-0.0) to (1.4,1);
\draw[fill = blue ,ultra thin] (0,0) circle(2pt) ;

\draw [->, bend left, thick]  (0.5,-0.0) to (1.5,-0.5);
\draw[fill = blue ,ultra thin] (0.5,0) circle(2pt) ;

\draw [bend left, thick, gray]   (-1,0) to (-1.4,0);

%\draw [bend right, thick]    (-1.4,0) to (0,0);
%\draw [bend right, thick]    (-1.4,0) to (0.5,0);
\draw [bend right, thick, gray]   (-1,0) to(0,0);
%\draw [bend right, thick]   (-1,0) to  (0.5,0);
\draw [bend right, thick, gray]  (0,0) to (0.5,0);

 \end{tikzpicture}
\caption{Realization of an \textit{open} street (thick red) with four neighboring streets. Four interconnected (connections drawn in grey) devices move towards the four neighboring streets (indicated by black arrows).}
\label{Fig_good_street}
\end{figure}
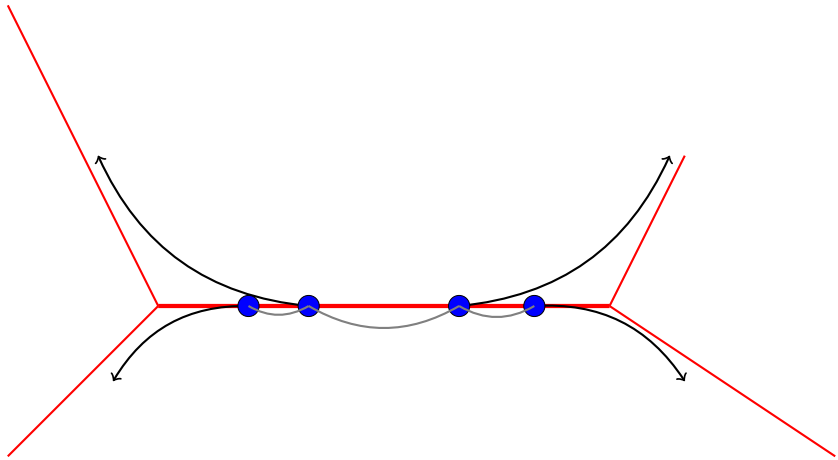

Let us note that we cannot assume that the probability for a street to be \textit{open} is monotone in the length of the street as there are for example the following two opposing effects:
If the length of the street increases, the probability of the event described in Condition (2) increases. However, the probability of the event described in Condition (1) is not monotone: While longer streets have more devices to establish connections, the distance between the crossings that needs to be bridged increases.
The property of a street to be open can be seen as a local event, since it depends only on the devices that start their movement on that street.
As $\k$ is well behaved, the probability to see a device that fulfills the role of one of the special device in Condition (2) is positive. Therefore, $\P(s \text{ is open})$ converges to one as $\lambda$ increases, by monotone convergence.

Now, if two adjacent streets $s$ and $s'$ are both open, due to Condition (2), there is both a node on $s$ that moves to street $s'$ and a node on street $s'$ that moves to street $s$. As they move in opposite directions, they will meet, and thus form a connection. The special devices on street $s$ are therefore connected to the special devices on street $s'$ and hence, if the street system that only consists of open streets contains an unbounded component, we obtain percolation of $g_{T,\muv,0,0}(X^\la)$.

Further, in order to establish percolation, in the second main step, we use a familiar discretization argument via large boxes in $\R^2$ in combination with stabilization and asymptotically essentially connectedness. 
Recall the random field of stabilization radii $\{R_x\}_{x\in\R^2}$ for the street system $S$. 
We say that $z \in \Z^2$ is \textit{$n$-open}, if we have that 
\begin{itemize}
\item[(a)] $R(Q_{6n}(nz))< n/2$,
\item[(b)] the number of streets that are fully contained in $Q_{6n}(nz)$ is smaller than $m$, and 
\item[(c)] every street fully contained in $Q_{6n}(nz)$ is open.
\end{itemize}
Now, if two adjacent boxes $Q_n(nz)$ and $Q_n(nw)$ are open, they are both non-empty due to the asymptotic essential connectedness, which is implied by Condition (a) since $R(Q_{6n}(nz))< n/2$ implies $R(Q_{2n}(nz))< n/2$. Additionally, under Condition (a), there is a path in $Q_{6n}(nz) \cap S$ connecting the streets in $Q_n(nz)$ and $Q_n(nw)$. In particular, due to Condition (c), all streets on this path are open. Therefore, the devices $Q_n(nz) \cap X^\la$ and $Q_n(nw)\cap X^\la$  are in the same connected component on $g_{T,\muv,0,0}(X^\la)$. Consequently, percolation of $n$-open boxes implies percolation of $g_{T,\muv,0,0}(X^\la)$.

In order to see percolation of $n$-open sites in $\Z^2$, note that the process of $n$-open sites, due to Condition (a) (which implies stabilization), is $6$-dependent. 
For $z=o$, let us denote by $A(n)$ the event described in Condition (a), by $B(m,n)$ the event described in Condition (b), and by $C(\la,m,n)$ the event described in Condition (c).
Thus, similar to the proof in Section\ref{Sec_Proof_Sub_2} we can invoke the domination-by-product-measures theorem~\cite[Theorem 0.0]{domProd} to establish Bernoulli site percolation on $\Z^2$ if 
\begin{align*}
\limsup_{n \uparrow\infty}\P( A^{\rm c}(n))=0 \qquad \text{ and for all $n\in \N$}\qquad\limsup_{m \uparrow\infty}\P( B^{\rm c}(m,n))=0
\end{align*}
and that for all $n,m\in \N$
\begin{align*}
\limsup_{\la \uparrow\infty}\P( C^{\rm c}(\la,m,n)\cap B(m,n))=0. 
\end{align*}
The first statement holds as $S$ is stabilizing. The second statement holds by the locally finiteness of the street system. Finally, for the third statement, note that, under the event $B(m,n)$, the number of fully contained streets is bounded by $m$ and hence we can apply the dominated convergence theorem together with the well-behavedness of the waypoint kernel to deduce
\begin{align*}
\limsup_{\la \uparrow\infty}\P( C^{\rm c}(\la,m,n)\cap B(m,n))\le \E\big[\one\{B(m,n)\}\sum_{s\in S\cap Q_{6n}}\limsup_{\la \uparrow\infty}\P(s\text{ is not open})\big]=0
\end{align*}
This finishes the proof. 
\end{proof}

\subsection{Proof of Theorem~\ref{thm-percolation-rho-positive}}\label{Sec_Proof_Sup_2}
\begin{proof}[Proof of Theorem~\ref{thm-percolation-rho-positive}]
The proof for this theorem is similar to the proof of Theorem~\ref{thm_largeintense}, although for $\rho > 0$ we have to add additional constraints to ensure that devices moving in opposite directions can actually connect as we can not control the moment devices turn around a corner, thus loosing connections.

To begin with, for a street $s$, let us denote by $s_1,\dots s_n$ the neighboring streets that share a joint crossing ($C_1, C_2$) with $s$.
Next, we will call $s$ \textit{open}, if there exist devices
$X^s_{i}, i \in \{1,\dots,n,n+1,n+2\}, $ on $s$ at $t = 0$ such that,
\begin{itemize}
\item[(1)] $\{X^s_{i}\}_{ i \in \{1,..,n+2\}}$ are connected in  $g_{T,\muv,\rho,r}(X^\la\cap s)$,
\item[(2)]  for all $i \in {1,\dots,n}$, there are times $t_i<\rho+\eps$ and $t_i< t'_i<T$ such that $\tra_{X_i^s,t} \in s$ for $t \in [0,t_i]$ and $\tra_{X_i^s,t} \in s_i$ for $t \in [t_i,t'_i]$ and
\item[(3)] For all $t \in [0,T]$ we have $\tra_{X_{n+1}^s,t} \in s$, $|\tra_{X_{n+1}^s,t} - C_1\
| < \rho v_{\rm min}$, $\tra_{X_{n+2}^s,t} \in s$, and $|\tra_{X_{n+2}^s,t} - C_2\
| < \rho v_{\rm min}$.
\end{itemize}
In words, Condition (1) ensures that the devices on $s$ are connected using only devices on $s$. The Condition (2) makes a statement for devices $X_i^s$, for $i \in 1,\dots ,n$ which play the role of 'transmitter' devices. They spend enough time on $s$ to become connected, by Condition (1), but not too much time, so that they can move to street $s_i$ and still spend enough time on $s_i$ to be able to 'transmit' there. The Condition (3) makes a statement about two devices that play the role of 'receiver' devices. They become connected, by Condition (1), and, due to their proximity to the crossing, are able to 'receive' a connection from an incoming device from a neighboring street.
%$s_i$ (if $s_i$ carries a 'transmitter' moving towards $s$, which is guaranteed for example if $s_i$ is open). 

Therefore, if two adjacent streets are both open, the 'transmitter' device that moves to the other street is able to connect to the 'receiver' device before time $T$ and vice versa. As a consequence, percolation of open streets implies percolation of $g_{T,\muv,\rho,r}(X^\la)$.

\medskip
Next, we use our conditions to ensure that streets are open with high probability as $\lambda$ increases and that open streets percolate. For this, first note that we assume that 
%For this, let us first assume that $v_{\rm min}<\min(r\rho^{-1}, a_{\rm c}(0)(2\rho)^{-1},c(2 \rho)^{-1})$, and note that,
%Recall that, by the definition of $v_{\rm min}$, there exists $\delta>0$, such that $\muv([v_{\rm min},v_{\rm min}+\delta'])>0$ for all $0<\delta'<\delta$. 
$2\rho v_{\rm min}<a_{\rm c}(0)$, and hence there exist $\delta>0$ and $0<\eps< T-2\rho$ such that also $2(\rho+\eps)(v_{\rm min}+\delta)<a_{\rm c}(0)$. This is useful since devices have to spend at least time $\rho+\epsilon$ on $s$ to satisfy Condition (2) and to connect with each other as demanded by Condition (1). This can only be satisfied for streets of length larger than $(\rho+\eps)(v_{\rm min}+\delta)$. However, a length of $2(\rho+\eps)(v_{\rm min}+\delta)$ is needed to ensure that the path of the device is a shortest path and moves through the joint crossing. 

Further note that for any pair of devices that moves in opposite directions on the same street, in order to be in contact for time $\rho$, the vicinity parameter needs to satisfy $r>\rho v_{\rm min}$, and this is the reason why we introduce this assumption. 

In particular, by possibly making $\delta$ even smaller, we have adjusted all the parameters $\rho>0$, $r>\rho (v_{\rm min}+\delta)$ and $T>2\rho$, such that the probability for any given street $s$ with $|s|>2(\rho+\eps)(v_{\rm min}+\delta)$ to be open tends to one as $\la$ tends to infinity. Here we also use the $c$-well-behavedness of the waypoint kernel, since this ensures that transmitting devices can be found. To see this, recall that the distance transmitting devices have to travel is at least $2\rho v_{\rm min}$. Also note, that the conditions are designed in such a way that, conditioned on $S$, openness of a given street $s\subset S$ can be checked independently of all the other streets in $S\setminus s$, and independently of the spatial position of the streets, by translation covariance of $\k$.

\medskip
Unfortunately, the probability of a street to be open depends on the whole surrounding street system, its length, its number of neighboring streets, its angle towards the neighboring streets and many more local properties of the street system. Therefore, there is no uniform bound for the probability that a street is open. 
However, we are able to guarantee two things: For each street $s$, with length larger than $2\rho v_{\rm min}$, we have $\lim_{\lambda \uparrow \infty} \P(s\text{ is open})=1$ and that streets with lengths larger than $2\rho v_{\rm min}$ percolate, by our assumption concerning $a_{\rm c}(0)$.
%While it is intuitively clear, that the process of open streets should percolate if $\lambda$ is chosen high enough, we are able to prove it with a reverse-Russo--Seymour--Welsh-type argument on a discretization:\\
Those arguments were sufficient in the proof of Theorem~\ref{thm_largeintense}, as we were able to use that the street system of open streets is asymptotically essentially connected. Unfortunately, the probability of streets with length smaller than $2\rho v_{\rm min}$ to be open may be zero. In order to prove percolation, we thus need a notion of asymptotically essentially connectedness on the thinned cluster.
This notion can be achieved either by definition, more precisely by the requirement $2\rho v_{\rm min}<a_{\rm c}^+$, or by choosing $v_{\rm min}$ small enough such that discrete arguments can use the asymptotic-essential connectedness of the underlying street system.
In the former case, the edges of length greater than $a_{\rm c}^+$ form an infinite cluster that is asymptotically essentially connected and the proof of Theorem~\ref{thm_largeintense} can be adapted to obtain percolation.

\medskip
In the following proof we will focus on the second approach of choosing $v_{\rm min}$ small enough, 
Let $(R_x)_{x\in \R^2}$ be the field of stabilization radii of the street system. We say that $z \in \Z^2$ is {\em $n$-open} if
\begin{itemize}
    \item[(a)] $R(Q_{3n}(nz)) < n$,
    \item[(b)] every street in $Q_{3n}(nz)$ has length larger than $2\rho v_{\rm min}$, and
    \item[(c)] every street in $S^{2\rho v_{\rm min}}\cap Q_{3n}(nz)$ is open.
\end{itemize}
Under those three conditions the existence of an infinite open path can be guaranteed by the original underlying street system.
Due to the $3$-dependence, guaranteed by Condition (a), we are able to use the domination-by-product-measure theorem~\cite[Theorem 0.0]{domProd} to establish Bernoulli site percolation on $\Z^2$, once we have shown that 
\begin{align*}
\limsup_{n\uparrow\infty}\limsup_{v_{\rm min}\downarrow0}\limsup_{\lambda\uparrow\infty}\P( o \text{ is not open} )=0.
\end{align*}
This is indeed the case. We label by $A(n), B(n, v_{\rm min}),$ and $C(n,v_{\rm min},\la)$ the events that the corresponding conditions are fulfilled for $z = o$. 
% As $\P(A_n^{\rm c}\cup B_n^{\rm c}\cup C_n^{\rm c})\le \P(A_n^{\rm c}\cup)...$
Then, note that 
$$\P( o \text{ is not open} )\le \P(A^{\rm c}(n))+\P(B^{\rm c}(n,v_{\rm min})\cap A(n))+\P(C^{\rm c}(n,v_{\rm min},\la)\cap B(n,v_{\rm min})\cap A(n)),$$
where $\limsup_{n\uparrow\infty}\P(A^{\rm c}(n))=0$ by stabilization. Further, for any fixed $n\in \N$, by Fatou's lemma,
$$\limsup_{v_{\rm min}\downarrow0}\P(B^{\rm c}(n,v_{\rm min})\cap A(n))\le\E\big[\limsup_{v_{\rm min}\downarrow0}\one\{B^{\rm c}(n,v_{\rm min})\cap A(n)\}\big] =0,$$
where we used that the street system is concentrated on streets of positive lengths. Finally, again by Fatou's lemma, for all $n\in \N$ and $v_{\rm min}>0$, we have that also
$$\limsup_{\lambda\uparrow\infty}\P(C^{\rm c}(n,v_{\rm min},\la)\cap B(n,v_{\rm min})\cap A(n)) =0.$$
Here, we also used that, under the event $B(n,v_{\rm min})$, every relevant street is sufficiently long and thus large $\la$ leads to openness of every street.
This finishes the proof. 
%Let us fix $\epsilon > 0$ arbitrarily small. Now, as the street system is stabilizing there is a $n_0(\epsilon)$ such that for all $n>n_0$ $\P(A_n^{\rm c})<\epsilon/3$. As the street system is locally finite, if $v_{\rm c}$ is chosen small enough  $\P(B_n^{\rm c})<\epsilon/3$. This choice is possible as the probability decreases monotonic in $v_{\rm min}$. Finally, there is a $\lambda_0=\lambda_0(\max(n_0,n_1),\epsilon)$ such that for $\lambda > \lambda_0$ we have $\P(C_n^{\rm c})<\epsilon/3$. Here, the probability to be open increases in $\lambda$. However, we are not able to give a uniform bound for $\lambda_{\rm c}$ as it depends on the previous choice of $v_{\rm min}$. As $\epsilon$ was chosen arbitrarily small, this completes the proof.
\end{proof}

\subsection{Proof of Theorem~\ref{thm_twicethinstreet}}\label{Sec_Proof_Sup_3}
\begin{proof}[Proof of Theorem~\ref{thm_twicethinstreet}]
We need to show that there exist $a > 0$  such that $S^{a}$ percolates with positive probability. 
Similarly to previous proofs, we construct a percolation process of finite-dependent open sites in $\Z^2$ that implies percolation of $S^{a}$.
For this, we say that $z \in \Z^2$ is \textit{$n$-open} if
\begin{itemize}
\item[(a)] $R(Q_{6n}(nz))< n/2$ and
\item[(b)] $Q_{6n}(nz)\cap S =  Q_{6n}(nz)\cap S^{a}$.
\end{itemize} 
Note that Condition (b) simply represents the event that no elimination of streets appears in $Q_{6n}(nz)$. 
Now the argument goes as follows. Due to Condition (a), by asymptotic essential connectedness, we have that $Q_n(nx) \cap S \neq \emptyset$ and that $Q_{3n}(nz)\cap S$ is connected within $Q_{6n}(nz)\cap S$. Then, for any two neighboring open sites $z,w\in \Z^2$, we have that $Q_n(nz)\cap S$ and  $Q_n(nw)\cap S$ are connected, as  $Q_n(nw) \subset Q_{3n}(nz)$. Next, due to Condition (b), all streets inside $Q_n(nz)$ and $Q_n(nw)$ lie in the same connected component of $S^{a}$. Therefore, percolation of open sites implies percolation of $S^{a}$.
Now, for the final argument we are able to use the domination-by-product-measures theorem~\cite[Theorem 0.0]{domProd} to establish supercritical Bernoulli site percolation on $\Z^2$ similar to the previous proofs.
% In order to establish the percolation regime for the process of $n$-good sites, note that, Condition (a)  additionally ensures stabilization and hence our process of $n$-good sites is $12$-dependent (and also translation invariant). Thus, similar to the previous proofs, we are allowed to use the domination-by-product-measures theorem~\cite[Theorem 0.0]{domProd} to establish Bernoulli site percolation on $\Z^2$ once we have shown that 
% \begin{align*}
% \limsup_{n,L\uparrow\infty}\limsup_{a,p \downarrow0}\P( o \text{ is not open} )=0.
% \end{align*}
% For this, let us denote by $A(n)$ the event described in Condition (a) and by $B(n,a,L,p)$ the event described in Condition (b). Then, it suffices to show that
% \begin{align*}
% \limsup_{n \uparrow\infty}\P( A^{\rm c}(n))=0,
% \end{align*}
% which is true as PVT is exponentially stabilizing, see~\cite[Example 3.1]{wias2}, and that for all $n\in \N$
% \begin{align*}
% \limsup_{L\uparrow\infty}\limsup_{a,p \downarrow0}\P( B^{\rm c}(n,a,L,p)\cap A(n))=0. 
% \end{align*}
% However, the thinning procedure is monotone in all parameters $a,L$ and $p$, and hence this statement since it is satisfied by weak convergence of $S^{a, L, p}$ towards $\emptyset$ as $a,p \downarrow0$ and $L\uparrow\infty$, and this finishes the argument.

\medskip
In order to prove that also $a_{\rm c}^+>0$, we are going to modify the definition of the random field of stabilizing radii $(R_x)_{x \in \R^2}$ and build a new stationary random field $(\bar R_x)_{x \in \R^2}$ based on the same probability space as $S$. For this, 
%let us denote by $C_\infty$ the infinite cluster of $S^a$
for any point $x \in \R^2$, recall that
\begin{align*}
    R_x^a &= \sup\{|y-x| \colon x\text{ and }y \text{ are connectable by a continuous path in $\R^2$ that does not intersect } S^a\},
%    &=\inf\{s> 0 \colon \text{there is a loop surrounding } x \text{ in } S^a\cap Q_s(x)\}
\end{align*}
and set $\bar R^a_x:= \max(R_x,R^a_x)$. We first prove that the field $\bar{R}^a_x$ is a legitimate alternative random field of stabilization radii in the sense of Definition~\ref{stabDef} and Definition~\ref{aecDef}. Since it is defined as a maximum, it suffices to show that $\lim_{n\uparrow\infty}\P(\bar R^a(Q_n) < n) = 1$ for some $a>0$, where again $\bar R^a(Q_n)=\sup_{x\in Q_n \cap\, \Q^2}\bar R^a_x$.  
%While the infinite cluster exists with probability one, we need global information on $S$ to discern where the infinite cluster ends. In contrast, $R^1$ only needs local information similar to the original radius $R$. Therefore, $(S,R^1)$ is jointly stationary and $R_1$ is measurable.

%As for all $x \in \Q$ we have $\bar R_x \ge R_x$ and we applied a local thinning rule, the independence of the original street system implies independence of the thinned street system. Similarly, as the original model was jointly stationary, so is the new one.
%Let us first prove the asymptotic-essentially connectedness. Due to stationarity it suffices to prove it for the origin. If $\sup_{y\in Q_n \cap \Q^2}\bar{R}_y \le n/2$, each point in $\Q_n$ is surrounded by a loop of the infinite cluster. The joint union of all these loops, is contained in $Q_{2n}$,as $R_y<1/2$. Per definition, there is no continuous path from $Q_n\cap\Q^2$ that leaves $Q_{2n}$ without crossing the infinite cluster. As the street system is locally finite and $Q_n\cap\Q^2$ is dense, this implies that there is no continuous path from $Q_n$ leaving $Q_{2n}$ without crossing the infinite cluster. Therefore, there exists a loop of $C_\infty$ in $Q_{2n}$ that surrounds $Q_n$.As the street system is planar, any street in $\Q_n$ that is part of the infinite cluster has to be connected to this surrounding loop. Therefore, $C_\infty\cap Q_n$ is connected in $C_\infty\cap Q_{2n}$.

%Now, it remains to prove that $\lim_{n\uparrow\infty}\P(\sup_{y\in Q_n \cap\, \Q^2}\bar R_y < n) = 1$.
We will prove this in two steps first using Bernoulli site percolation on $\Z^2$ and then a Russo--Seymour--Welsh-type argument. 
%statement to show that the probability of a left-right crossing goes to one.
Let us fix $n>0$.
We say that $z \in \Z^2$ is {\em $n$-good} if
\begin{itemize}
     \item[(a)] $R(Q_n) < n$ and
     \item[(b)] $S\cap Q_{6n}(nz) = S^a \cap Q_{6n}(nz)$.
\end{itemize}
As $S$ is asymptotically essentially connected and Condition (b) guarantees that no streets were removed, adjacent $n$-good sites guarantee that the street system in the associated $n$-squares $Q_n$ are connected. 
%Therefore, an infinite cluster of open boxes implies that neighboring open boxes are able to connect to each other.
As the good sites are $6$-dependent, we can link it with the domination-by-product-measures theorem~\cite[Theorem 0.0]{domProd} to establish Bernoulli site percolation on $\Z^2$ once we have shown that $\limsup_{n \uparrow \infty}\limsup_{a \downarrow 0}\P(o \text{ is not } n\text{-good})= 0$.
This is indeed the case, as we can first take $n$ large enough such that Condition (a) is satisfied with high probability, and later choose $a(n)$ small enough such that the probability of the joint event is arbitrarily close to one.
Now, as we consider Bernoulli site percolation, we can apply a Russo--Seymour--Welsh-type argument as follows.
The probability to see a left-right crossing in the lattice of open vertices in a $m \times 3m$ box converges to one for $m\uparrow \infty$.
Therefore, by surrounding a box of size $m$ with four of those corridors from each side, we see that with high probability, the $m$-box is surrounded, see Figure~\ref{figRSW}.
\begin{figure}\label{figRSW}
\centering
\begin{tikzpicture}[scale = 2.0]
\draw[thick] (0.5,-1.5) -- (0.5,1.5);
\draw[thick] (-0.5,-1.5) -- (-0.5,1.5);
\draw[thick] (1.5,-1.5) -- (1.5,1.5);
\draw[thick] (-1.5,-1.5) -- (-1.5,1.5);
\draw[thick] (-1.5,1.5) -- (+1.5,1.5);
\draw[thick] (-1.5,0.5) -- (+1.5,0.5);
\draw[thick] (-1.5,-0.5) -- (+1.5,-0.5);
\draw[thick] (-1.5,-1.5) -- (+1.5,-1.5);
\draw [red, very thick] plot [smooth, tension=1] coordinates {(1.5,1.2) (1.2,1.3) (1.1,1.3) (1.0,0.9) (0.8,1.1) (0.5,1.2) (0.1,0.6) (-0.5,1.2) (-1.1,1.1) (-1.5,0.9)};
\draw [red, very thick] plot [smooth, tension=1] coordinates{(1.1,-1.5) (0.9,-1.1) (1.3,-0.9) (1.1, -0.6) (0.7, -0.2) (1.0,0.4) (1.1,1.0) (1.3,1.2)(1.2,1.5)};
\draw [red, very thick] plot [smooth, tension=1] coordinates {(1.5,-1.2) (1.2,-1.3) (1.0,-0.9) (0.8,-1.1) (0.5,-1.1) (0.1,-0.7) (-0.5,-0.7) (-1.1,-1.2) (-1.5,-0.6)};
\draw [red, very thick] plot [smooth, tension=0.5] coordinates{(-0.7,1.5)(-1.2,1.3)(-1.2,1.0)(-1.4,0.6)(-0.7,0.2)(-1,-0.2)(-1.3,-0.6)(-0.6,-1.5)};
%\draw[ultra thick, blue] (-1.5,1.5)--(0.5,-0.5);

\draw[ultra thick, fill = black!20] (-0.5,-0.5) -- (0.5,-0.5) -- (0.5,0.5) -- (-0.5,0.5) -- cycle ;
\draw[decorate, decoration={brace}, yshift=0.5ex]  (-1.5,1.5) -- node[above=0.4ex] {$3mn$}  (1.5,1.5);
\end{tikzpicture}
\caption{Illustration of a good box: The gray box consists of $m^2$ boxes $[0,n]^2$. It is surrounded by left-right and up-down crossings of $n$-good boxes.}
\end{figure}
Now, that the $m$-box is surrounded, we have a joint bound for $\bar R^a_x<2\sqrt{2}mn$ for all $x \in Q_{mn}$. Putting things together, we first pick $n$ sufficiently large and $a$ sufficiently small, such that $\lim_{m\uparrow\infty}\P(\bar R^a(Q_{mn})\ge 2\sqrt{2}mn)=0$. Then, using a union bound and shift invariance, this implies that also $\lim_{m\uparrow\infty}\P(\bar R^a(Q_{mn})\ge mn)=0$ as required.
\end{proof}

%====================================================
%====================================================

\bibliography{wias}
\bibliographystyle{alpha}

\end{document}